\documentclass[namedreferences,11pt,letterpaper]{kluwer}

\usepackage{amsfonts}
\usepackage{times}
\usepackage{mathptm}
\usepackage{german}
\usepackage{kluref}
\usepackage{proof}
\def\gq#1{\selectlanguage{german}#1\selectlanguage{USenglish}}
\def\sg{\selectlanguage{german}}
\def\eg{\selectlanguage{USenglish}}
\let\citeA\opencite
\let\citeN\inlinecite
\def\pciteN#1#2{\citeauthor{#1}~(\citeyear{#1}, #2)}

\def\pcite#1#2{(\citeauthor{#1}~\citeyear{#1}, #2)}

\let\note\endnote
\let\ol\overline
\def\Bcite#1#2{\pciteN{#2}{#1}}
\def\pp#1{#1}
\let\epsilon\varepsilon
\let\phi\varphi
\def\tfA{\tilde{\frak{A}}}
\def\tfB{\tilde{\frak{B}}}
\def\fz{\frak{z}}

\begin{document}                 

\addtolength{\topmargin}{2cm}                                                                  
\begin{article}
\begin{opening}         
\selectlanguage{USenglish}

\title{The Practice of Finitism: 
Epsilon Calculus and Consistency Proofs in Hilbert's Program}
\author{Richard \surname{Zach}\email{zach@math.berkeley.edu}}
\institute{University of California, Berkeley}
\runningtitle{The Practice of Finitism}
\date{Second Draft, February 21, 2001 -- Comments welcome!}

\begin{ao}
Group in Logic and the Methodology of Science\\
University of California, Berkeley\\
Berkeley, CA 94720--3840\\
http://www.math.berkeley.edu/$\sim$zach/
\end{ao}

\begin{abstract}
After a brief flirtation with logicism in 1917--1920, David Hilbert
proposed his own program in the foundations of mathematics in 1920 and
developed it, in concert with collaborators such as Paul Bernays and
Wilhelm Ackermann, throughout the 1920s.  The two technical pillars of
the project were the development of axiomatic systems for ever
stronger and more comprehensive areas of mathematics and finitistic
proofs of consistency of these systems.  Early advances in these areas
were made by Hilbert (and Bernays) in a series of lecture courses at
the University of G\"ottingen between 1917 and 1923, and notably in
Ackermann's dissertation of 1924.  The main innovation was the
invention of the $\epsilon$-calculus, on which Hilbert's axiom systems
were based, and the development of the $\epsilon$-substitution method
as a basis for consistency proofs.  The paper traces the development
of the ``simultaneous development of logic and mathematics'' through
the $\epsilon$-notation and provides an analysis of Ackermann's
consistency proofs for primitive recursive arithmetic and for the
first comprehensive mathematical system, the latter using the
substitution method.  It is striking that these proofs use transfinite
induction not dissimilar to that used in Gentzen's later consistency
proof as well as non-primitive recursive definitions, and that these
methods were accepted as finitistic at the time.
\end{abstract}

\end{opening}


\bibliographystyle{klunamed}


\section{Introduction}\label{sec-conprf-intro}

Hilbert first presented his philosophical ideas based on the axiomatic
method and consistency proofs in the years 1904 and 1905, following
his exchange with Frege on the nature of axiomatic systems and the
publication of Russell's Paradox. In the text of Hilbert's address to the
International Congress of Mathematicians in Heidelberg, we read:
\begin{quotation}
Arithmetic is often considered to be part of logic, and the
traditional fundamental logical notions are usually presupposed when
it is a question of establishing a foundation of arithmetic.  If we
observe attentively, however, we realize that in the traditional
exposition of the laws of logic certain fundamental arithmetic notions
are already used, for example, the notion of set and, to some extent,
also that of number. Thus we find ourselves turning in a circle, and
that is why a partly simultaneous development of the laws of logic and
arithmetic is required if paradoxes are to be
avoided.\note{\pciteN{Hilbert:05}{131}.  For a general discussion of
Hilbert's views around 1905, see \pciteN{Peckhaus:90}{Chapter~3}.}
\end{quotation}
When Hilbert returned to his foundational work with full force in 1917,
he seems at first to have been impressed with Russell's and Whiteheads's
work in the \emph{Principia}, which succeeded in developing large parts 
of mathematics without using sets.  By 1920, however, Hilbert returned
to his earlier conviction that a reduction of mathematics to logic is
not likely to succeed.  Instead, he takes Zermelo's axiomatic set theory
as a suitable framework for developing mathematics.  He localizes the
failure of Russell's logicism in its inability to provide the existence 
results necessary for analysis:
\begin{quotation}
The axiomatic method used by Zermelo is unimpeachable and
indispensable.  The question whether the axioms include a
contradiction, however, remains open.  Furthermore the question poses
itself if and in how far this axiom system can be deduced from logic.
[\dots T]he attempt to reduce set theory to logic seems promising
because sets, which are the objects of Zermelo's axiomatics, are
closely related to the predicates of logic.  Specifically, sets can be
reduced to predicates.

This idea is the starting point for Frege's, Russell's,
and Weyl's investigations into the foundations of mathematics.\note{
\sg"`Die von {\sc zermelo} benutze axiomatische Methode ist zwar
unanfechtbar und unentbehrlich.  Es bleibt dabei doch die Frage offen,
ob die aufgestellten Axiome nicht etwa einen Widerspruch
einschliessen.  Ferner erhebt sich die Frage, ob und inwieweit sich
das Axiomensystem aus der Logik ableiten l"asst. [\dots D]er Versuch
einer Zur"uckf"uhrung auf die Logik scheint besonders deshalb
aussichtsvoll, weil zwischen Mengen, welche ja die Gegenst"ande in
{\sc zermelo}s Axiomatik bilden, und den Pr"adikaten der Logik ein
enger Zusammenhang besteht.  N"amlich die Mengen lassen sich auf
Pr"adikate zur"uckf"uhren.

Diesen Gedanken haben {\sc frege}, {\sc russel[l]} und {\sc weyl} zum
Ausgangspunkt genommen bei ihren Untersuchungen "uber die Grundlagen
der Mathematik."' \pcite{Hilbert:20b}{27--28}.}
\end{quotation}
The logicist project runs into a difficulty when, given a second-order
predicate $S$ to which a set of sets is reduced, we want to know that
there is a predicate to which the union of the sets reduces.  This
predicate would be $(\exists P)(P(x) \& S(P))$---$x$ is in the union
of the sets in $S$ if there is a set $P$ of which $x$ is a member and
which is a member of $S$.
\begin{quotation}
We have to ask ourselves, what ``there is a predicate $P$'' is
supposed to mean.  In axiomatic set theory ``there is'' always refers
to a basic domain $\frak{B}$.  In logic we could also think of the
predicates comprising a domain, but this domain of predicates cannot
be seen as something given at the outset, but the predicates must be
formed through logical operations, and the rules of construction
determine the domain of predicates only afterwards.

From this we see that in the rules of logical construction of
predicates reference to the domain of predicates cannot be allowed.
For otherwise a \emph{circulus vitiosus} would result.\note{\sg"`Wir
m"ussen uns n"amlich fragen, was es bedeuten soll: "`es gibt ein
Pr"adikat $P$."' In der axiomatischen Mengentheorie bezieht sich das
"`es gibt"' immer auf den zugrunde gelegten Bereich $\frak{B}$.  In
der Logik k"onnen wir zwar auch die Pr"adikate zu einem Bereich
zusammengefasst denken; aber dieser Bereich der Pr"adikate kann hier
nicht als etwas von vorneherein Gegebenes betrachtet werden, sondern
die Pr"adikate m"ussen durch logische Operationen gebildet werden, und
durch die Regeln der Konstruktion bestimmt sich erst nachtr"aglich der
Pr"adikaten-Bereich.

Hiernach ist ersichtlich, dass bei den Regeln der logischen
Konstruktion von Pr"adikaten die Bezugnahme auf den
Pr"adikaten-Bereich nicht zugelassen werden kann. Denn sonst erg"abe
sich ein \emph{circulus vitiosus}."' \pcite{Hilbert:20b}{31}.}
\end{quotation}
Here Hilbert is echoing the predicativist worries of Poincar\'e and
Weyl.  However, Hilbert rejects Weyl's answer to the problem, viz.,
restricting mathematics to predicatively acceptable constructions and
inferences, as unacceptable in that it amounts to ``a return to the
prohibition policies of Kronecker.''  Russell's proposed solution, on
the other hand, amounts to giving up the aim of reduction to logic:
\begin{quotation}
Russell starts with the idea that it suffices to replace the
predicate needed for the definition of the union set by one
that is extensionally equivalent, and which is not open to
the same objections.  He is unable, however, to exhibit
such a predicate, but sees it as obvious that such a predicate
exists.  It is in this sense that he postulates the ``axiom
of reducibility,'' which states approximately the following:
``For each predicate, which is formed by referring (once or multiple times)
to the domain of predicates, there is an extensionally
equivalent predicate, which does not make such reference.

With this, however, Russell returns from constructive logic to the
axiomatic standpoint. [\dots]

The aim of reducing set theory, and with it the usual methods of
analysis, to logic, has not been achieved today and maybe cannot be
achieved at all.\note{``\sg{\sc russell} geht von dem Gedanken aus,
dass es gen"ugt, das zur Definition der Vereinigungsmenge unbrauchbare
Pr"adikat durch ein sachlich gleichbedeutendes zu ersetzen, welches
nicht dem gleichen Einwande unterliegt.  Allerdings vermag er ein
solches Pr"adikat nicht anzugeben, aber er sieht es als ausgemacht an,
dass ein solches existiert.  In diesem Sinne stellt er sein "`Axiom
der Reduzierbarkeit"' auf, welches ungef"ahr folgendes besagt: "`Zu
jedem Pr"adikat, welches durch (ein- oder mehrmalige) Bezugnahme auf
den Pr"adikatenbereich gebildet ist, gibt es ein sachlich
gleichbedeutendes Pr"adikat, welches keine solche Bezugnahme aufweist.

Hiermit kehrt aber {\sc russell} von der konstruktiven Logik zu dem
axiomatischen Standpunkt zur"uck. [\dots]

Das Ziel, die Mengenlehre und damit die
gebr"auchlichen Methoden der Analysis auf die Logik zur"uckzuf"uhren,
ist heute nicht erreicht und ist vielleicht "uberhaupt nicht
erreichbar.\eg''
\pcite{Hilbert:20b}{32--33}.}
\end{quotation}
With this, Hilbert rejects the logicist position as failed.  At the
same time, he rejects the restrictive positions of Brower, Weyl, and
Kronecker.  The axiomatic method provides a framework which can
accommodate the positive contributions of Brouwer and Weyl, without
destroying mathematics through a Kroneckerian ``politics of
prohibitions.'' For Hilbert, the unfettered progress of mathematics,
and science in general, is a prime concern.  This is a position that
Hilbert had already stressed in his lectures before the 1900 and 1904
International Congresses of Mathematics, and which is again of
paramount importance for him with the conversion of Weyl to Brouwer's
intuitionism.

Naturally, the greater freedom comes with a price attached: the
axiomatic method, in contrast to a foundation based on logical
principles alone, does not itself guarantee consistency.  Thus, a
proof of consistency is needed.   


\section{Early Consistency Proofs}\label{sec-early-conprfs}

Ever since his work on geometry in the 1890s, Hilbert had an interest
in consistency proofs.  The approaches he used prior to the
foundational program of the 1920s were almost always relative
consistency proofs.  Various axiomatic systems, from geometry to
physics, were shown to be consistent by giving arithmetical (in a
broad sense, including arithmetic of the reals) interpretations for
these systems, with one exception---a prototype of a finitistic
consistency proof for a weak arithmetical system in
\citeN{Hilbert:05}.  This was Hilbert's first attempt at a ``direct''
consistency proof for arithmetic, i.e., one not based on a reduction
to another system, which he had posed as the second of his famous list
of problems \cite{Hilbert:00a}.

When Hilbert once again started working on foundational issues
following the war, the first order of business was a formulation of
logic.  This was accomplished in collaboration with Bernays between
1917 and 1920 (see \citeA{Sieg:99} and \citeA{Zach:99}), included the
establishment of metatheoretical results like completeness,
decidability, and consistency for propositional logic in 1917/18, and
was followed by ever more nuanced axiom systems for propositional and
predicate logic.  This first work in purely logical axiomatics was
soon extended to include mathematics.  Here Hilbert followed his
own proposal, made first in 1905,\note{``\dots zur Vermeidung von
Paradoxien ist daher eine teilweise gleichzeitige Entwicklung der
Gesetze der Logik und der Arithmetik erforderlich.''
\pcite{Hilbert:05}{176}.} to develop mathematics and logic
simultaneously.  The extent of this simultaneous development is
nowhere clearer than in Hilbert's lecture course of 1921/22, where the
$\epsilon$-operator is first used as both a logical notion,
representing the quantifiers, and an arithmetical notion, representing
induction in the form of the least number principle. Hilbert realized
then that a consistency proof for all of mathematics is a difficult
undertaking, best attempted in stages:
\begin{quotation}
Considering the great variety of connectives and interdependencies
exhibited by arithmetic, it is obvious from the start that we will
the not be able to solve the problem of proving consistency in one
fell swoop.  We will instead first consider the simplest connectives,
and then proceed to ever higher operations and inference methods,
whereby consistency has to be established for each extension of the
system of signs and inference rules, so that these extensions do not
endanger the consistency [result] established in the preceding stage.

Another important aspect is that, following our plan for the complete
formalization of arithmetic, we have to develop the proper
mathematical formalism in connection with the formalism of the logical
operations, so that---as I have expressed it---a simultaneous
construction of mathematics and logic is executed.\note{``\sg In
Anbetracht der grossen Mannigfaltigkeit von Verkn"upfungen und
Zusammenh"angen, welche die Arithmetik aufweist, ist es von vornherein
ersichtlich, dass wir die Aufgabe des Nachweises der
Widerspruchslosigkeit nicht mit einem Schlage l"osen k"onnen.  Wir
werden vielmehr so vorgehen, dass wir zun"achst nur die einfachsten
Verkn"upfungen betrachten und dann schrittweise immer h"ohere
Operationen und Schlussweisen hinzunehmen, wobei dann f"ur jede
Erweiterung des Systems der Zeichen und der Uebergangsformeln einzeln
der Nachweis zu erbringen ist, dass sie die auf der vorherigen Stufe
festgestellte Widerspruchsfreiheit nicht aufheben.

Ein weiterer wesentlicher Gesichtspunkt ist, dass wir, gem"ass unserem
Plan der restlosen Formalisierung der Arithmetik, den eigentlich
mathematischen Formalismus im Zusammenhang mit dem Formalismus der
logischen Operationen entwickeln m"ussen, sodass---wie ich es
ausgedr"uckt habe---ein simultaner Aufbau von Mathematik und Logik
ausgef"uhrt wird.\eg'', \pcite{Hilbert:21}{8a--9a}.  The passage is
not contained in Kneser's notes \cite{Hilbert:22c} to the same
course.}
\end{quotation}
Hilbert had rather clear ideas, once the basic tools both of proof and
of formalization were in place, of what the stages should be.  In an
addendum to the lecture course on {\em Grundlagen der Mathematik},
taught by Hilbert and Bernays in 1922--23,\note{The notes by Kneser
\cite{Hilbert:22c} do not contain the list of systems below.  The
version of the $\epsilon$-calculus used in the addendum is the same as
that used in Kneser's notes, and differs from the presentation in
\citeN{Ackermann:24}, submitted February 20, 1924.}  he outlined them.
The first stage had already been accomplished: Hilbert gave
consistency proofs for calculi of propositional logic in his 1917/18
lectures.  Stage~II consist in the elementary calculus of free
variables, plus equality axioms and axioms for successor and
predecessor.  The axioms are:
\[
\begin{array}{lclc}
1. & A \to B \to A & 
2. & (A \to A \to B) \to A \to B \\
3. & (A \to B \to C) \to (B \to A \to C) & 
4. & (B \to C) \to (A \to B) \to A \to C\\
5. & A \& B \to A &
6. & A \& B \to B \\
7. & A \to B \to A \& B &
8. & A \to A \lor B \\
9. & B \to A \lor B &
10. & (A \to C) \to (B \to C) \to A \lor B \to C \\
11. & A \to \overline{A} \to B &
12. & (A \to B) \to (\overline{A} \to B) \to B \\
13. & a = a &
14. & a = b \to A(a) \to A(b) \\
15. & a + 1 \neq 0 &
16. & \delta(a+1) = a\note{\cite{Hilbert:22b}, 17, 19 }
\end{array}
\] 
In Hilbert's systems, Latin letters are variables; in particular, $a$,
$b$, $c$, \dots, are individual variables and $A$, $B$, $C$, \dots,
are formula variables.  The rules of inference are modus ponens and
substitution for individual and formula variables.

Hilbert envisaged his foundational project as a stepwise
``simultaneous development of logic and mathematics,'' in which
axiomatic systems for logic, arithmetic, analysis, and finally set
theory would be developed.  Each stage would require a proof of
consistency before the next stage is developed.  In a handwritten
supplement to the typescript of the 1922--23 lecture notes on the
foundations of arithmetic, Hilbert presents a rough overview of what
these steps might be:
\begin{quotation}
{\bf Outline}.  Stage II was elementary calculation, axioms 1--16.

Stage III. Now elementary number theory\\
Schema for definition of functions by recursion and modus ponens\\
will add the schema of induction to modus ponens\\
even if this coincides in substance with the results of 
intuitively obtained number theory, we are now dealing with formulas, e.g, $a + b = b + a$.

Stage IIII.  Transfinite inferences and parts of analysis

Stage V. Higher-order variables and set theory. Axiom of choice.

Stage VI. Numbers of the 2nd number class, full transfinite induction.
Higher types. Continuum problem, transfinite induction for numbers in
the 2nd number class.

Stage VII. (1) Replacement of infinitely many definitional schemata by one axiom.  (2)
Analysis and set theory.  At level 4, again the full theorem of the least upper bound.

Stage VIII. Formalization of well ordering.\note{\sg``{\bf Disposition}.
Stufe II war elementares Rechnen Axiome 1--16

Stufe III. Nun elementare Zahlentheorie\\ Schema f"ur Def. von
Funktionen durch Rekursion u. Schlussschema\\{} wollen [?] unser
Schlussschema noch das Induktionsschema hinzuziehen\\ Wenn auch
inhaltlich das wesentlich mit den Ergebnissen der anschauliche
gewonnenen [?] Zahlenth.{} "ubereinstimmt, so doch jetzt Formeln
z.B. $a + b = b + a$.

Stufe IIII. Transfinite Schlussweise u. teilweise Analysis

Stufe V. H"ohere Variablen-Gattungen u. Mengenlehre. Auswahlaxiom

Stufe VI. Zahlen d[er] 2$^{\hbox{\footnotesize ten}}$ Zahlkl[asse], Volle transfin[ite] Induktion. H"ohere Typen. Continuumsproblem, transfin[ite] Induktion f"ur Zahlen der 2$^{\hbox{\footnotesize ten}}$ Zahlkl.

Stufe VII. 1.) Ersetzung der $\infty$ vielen Definitionsschemata durch ein Axiom.
2.) Analysis u[nd] Mengenlehre. Auf der 4$^{\rm ten}$ Stufe nochmals der volle Satz von der oberen Grenze

Stufe VIII. Formalisierung der Wohlordnung.
\eg'' \cite{Hilbert:22c}, \emph{Erg"anzung}, sheet 1.}
\end{quotation}

\subsection{The propositional calculus and the calculus of 
elementary computation}

Step I had been achieved in 1917--18.  Already in the lectures from
the Winter term 1917/18, Hilbert and Bernays had proved that the
propositional calculus is consistent.  This was done first by
providing an arithmetical interpretation, where they stressed that
only finitely many numbers had to be used as ``values'' (0 and 1).
The proof is essentially a modern proof of the soundness of
propositional logic: A truth value semantics is introduced by
associating which each formula of the propositional calculus a truth
function mapping tuples of 0 and 1 (the values of the propositional
variables) to 0 or 1 (the truth value of the formula under the
corresponing valuation).  A formula is called {\em correct} if it
corresponds to a truth function which always takes the value~1. It is
then showed that the axioms are correct, and that modus ponens
preserves correctness.  So every formula derivable in the
propositional calculus is correct.  Since $A$ and $\ol A$ cannot both
be correct, they cannot both be derivable, and so the propositional
calculus is consistent.

It was very important for Hilbert that the model for the propositional
calculus thus provided by $\{0, 1\}$ was finite.  As such, its
existence, and the admissibility of the consistency proof was beyond
question.  This lead him to consider the consistency proof for the
propositional calculus to be the prime example for for a consistency
proof {\em by exhibition} in his 1921/22 lectures on the foundations
of mathematics.  The consistency problem in the form of a demand for a
consistency proof for an axiomatic system which neither proceeds by
exhibiting a model, nor by reducing consistency of a system to the
consistency of another, but by providing a metamathematical proof that
no derivation of a contradiction is possible, is first formulated in
lectures in the Summer term of 1920.  Here we find a first formulation
of an arithmetical system and a proof of consistency.  The system
consists of the axioms 
\begin{eqnarray*} 
1 & = & 1\\ 
(a = b) & \to & (a + 1 = b + 1)\\ 
(a + 1 = b + 1) & \to & (a = b)\\ 
(a = b) & \to & ((a = c) \to (b = c)).
\end{eqnarray*}
The notes contain a proof that these four axioms, together with modus
ponens, do not allow the derivation of the formula
\[
a + 1 = 1.
\]
The proof itself is not too interesting, and I will not reproduce it
here.\note{The proof can also be found in
\pciteN{Hilbert:22a}{171--173}; cf.~\pciteN{Mancosu:98}{208--210}.}
The system considered is quite weak. It does not even contain all of
propositional logic: negation only appears as inequality, and only
formulas with at most two `$\to$' signs are derivable.  Not even $a=a$
is derivable.  It is here, nevertheless, that we find the first
statement of the most important ingredient of Hilbert's project,
namely, proof theory:
\begin{quotation}
Thus we are led to make the proofs themselves the object of our
investigation; we are urged towards a {\em proof theory}, which
operates with the proofs themselves as objects.

For the way of thinking of ordinary number theory the numbers are then
objectively exhibitable, and the proofs about the numbers already
belong to the area of thought.  In our study, the proof itself is
something which can be exhibited, and by thinking about the proof we
arrive at the solution of our problem.

Just as the physicist examines his apparatus, the astronomer his
position, just as the philosopher engages in critique of reason, so
the mathematician needs his proof theory, in order to secure each
mathematical theorem by proof critique.\note{``\sg Somit sehen wir uns
veranlasst, die Beweise als solche zum Gegenstand der Untersuchung zu
machen; wir werden zu einer Art von \emph{Beweistheorie} gedr"angt,
welche mit den Beweisen selbst als Gegenst"anden operiert.

F"ur die Denkweise der gew"ohnlichen Zahlentheorie sind die Zahlen das
gegenst"andlich-Aufweisbare, und die Beweise der S"atze "uber die
Zahlen fallen schon in das gedankliche Gebiet. Bei unserer
Untersuchung ist der Beweis sebst etwas Aufweisbares, und durch das
Denken "uber den Beweis kommen wir zur L"osung unseres Problems.

Wie der Physiker seinen Apparat, der Astronom seinen Standort
untersuch, wie der Philosoph Vernunft-Kritik "ubt, so braucht der
Mathematiker diese Beweistheorie, um jeden mathematischen Satz durch
eine Beweis-Kritik sicherstellen zu k"onnen.\eg''
\pciteN{Hilbert:20b}{39--40}.  Almost the same passage is found in
\pciteN{Hilbert:22a}{169--170}, cf.~\pciteN{Mancosu:98}{208}.}
\end{quotation}

This project is developed in earnest in two more lecture courses in
1921--22 and 1922--23.  These lectures are important in two respects.
First, it is here that the axiomatic systems whose consistency is to
be proven are developed.  This is of particular interest for an
understanding of the relationship of Hilbert to Russell's project in
the \emph{Principia} and the influence of Russell's work both on
Hilbert's philosophy and on the development of axiomatic systems for
mathematics.\note{For a detailed discussion of these influences, see
\citeN{Mancosu:99}.}  \citeN{Sieg:99} has argued that, in fact,
Hilbert was a logicist for a brief period around the time of his paper
``Axiomatic Thought'' \cite{Hilbert:18}.  However, as noted in
Section~\ref{sec-conprf-intro}, Hilbert soon became critical of
Russell's type theory, in particular of the axiom of reducibility.
Instead of taking the system of \emph{Principia} as the adequate
formalization of mathematics the consistency of which was to be shown,
Hilbert proposed a new system.  The guiding principle of this system
was the ``simultaneous development of logic and mathematics''---as
opposed to a development of mathematics out of logic---which he had
already proposed in \pciteN{Hilbert:05}{176}.  The cornerstone of
this development is the $\epsilon$-calculus.  The second major
contribution of the 1921--22 and 1922--23 lectures are the consistency
proofs themselves, including the \emph{Hilbertsche Ansatz} for the
$\epsilon$-substitution method, which were the direct precursors to
Ackermann's dissertation of 1924.

In contrast to the first systems of 1920, here Hilbert uses a system
based on full propositional logic with axioms for equality, i.e., the
elementary calculus of free variables:
\begin{quotation}
\noindent I. Logical axioms

\noindent a) Axioms of consequence

\noindent 1) $A \to B \to A$\\
2) $(A \to A \to B) \to A \to B$\\
3) $(A \to B \to C) \to B \to A \to C$\\
4) $(B \to C) \to (A \to B) \to A \to C$\\

\noindent b) Axioms of negation

\noindent 5) $A \to \ol{A} \to B$\\
6) $(A \to B) \to (\ol{A} \to B) \to B$

\noindent II. Arithmetical axioms

\noindent a) Axioms of equality

\noindent 7) $a=a$\\
8) $a = b \to Aa \to Ab$

\noindent b) Axioms of number

\noindent 9) $a + 1 \neq 0$\\ 10) $\delta(a + 1) =
a$\note{\pcite{Hilbert:21}{part 2, 3}.  Kneser's \emph{Mitschrift} of
these lectures contains a different system which does not include
negation.  Instead, numerical inequality is a primitive. This system
is also found in Hilbert's first talks on the subject in Copenhagen
and Hamburg in Spring and Summer of 1921.  \citeN{Hilbert:23}, a talk
given in September 1922, and Kneser's notes to the course of Winter
Semester 1922--23 \cite{Hilbert:22b} do contain the new system with
negation.  This suggests that the developments of Hilbert's 1921--22
lectures were not incorporated into the published version of Hilbert's
Hamburg talk~\shortcite{Hilbert:22a}. Although \shortcite{Hilbert:22a}
was published in 1922, and a footnote to the title says ``This
communication is essentially the content of the talks which I have
given in the Spring of this year in Copenhagen [\dots] and in the
Summer in Hamburg [\dots],'' it is clear that the year in question is
1921, when Hilbert addressed the Mathematisches Seminar of the
University of Hamburg, July 25--27, 1921.  A report of the talks was
published by Reidemeister in \emph{Jahrbuch der Deutschen
Mathematiker-Vereinigung}~\textbf{30}, 2.~Abt. (1921), 106.
\citeN{Hilbert:22b} also have separate axioms for conjunction and
disjunction, while in \shortcite{Hilbert:23} it is extended it by
quantifiers.}
\end{quotation}
Here, `${}+1$' is a unary function symbol.  The rules of inference are
substitution (for individual and formula variables) and modus ponens.

Hilbert's idea for how a finitistic consistency proof should be
carried out is first presented here.  The idea is this: suppose a
proof of a contradiction is available.  We may assume that the end
formula of this proof is $0 \neq 0$. 
\begin{enumerate}
\item \emph{Resolution into proof threads.} First, we observe that by
duplicating part of the proof and leaving out steps, we can transform
the derivation to one where each formula (except the end formula) is
used exactly once as the premise of an inference.  Hence, the proof is
in tree form.  
\item \emph{Elimination of variables.} We transform the proof so that
it contains no free variables.  This is accomplished by proceeding
backwards from the end formula: The end formula contains no free
variables.  If a formula is the conclusion of a substitution rule, the
inference is removed.  If a formula is the conclusion of modus ponens
it is of the form
\[
\infer{\frak{B'}}{\frak{A} & \frak{A} \to \frak{B}}
\]
where $\frak{B}'$ results from $\frak{B}$ by substituting terms for
free variables.  If these variables also occur in $\frak{A}$, we
substitute the same terms for them.  Variables in $\frak{A}$ which do
not occur in $\frak{B}$ are replaced with~$0$.  This yields a formula
$\frak{A}'$ not containing variables.\note{The procedure whereby we
pass from $\frak{A}$ to $\frak{A'}$ is simple in this case, provided
we keep track of which variables are substituted for below the
inference.  In general, the problem of deciding whether a formula is a
substitution instance of another, and to calculate the substitution
which would make the latter syntactically identical to the former is
known as \emph{matching}.  Although not computationally difficult, it
is not entirely trivial either.}  The inference is replaced by
\[
\infer{\frak{B'}}{\frak{A}' & \frak{A}' \to \frak{B}'}
\]
\item \emph{Reduction of functionals.}  The remaining derivation
contains a number of terms (\emph{functionals} in Hilbert's parlance)
which now have to be reduced to numerical terms (i.e., standard
numerals of the form $(\ldots(0 + 1) + \cdots) + 1$).  In this case,
this is done easily by rewriting innermost subterms of the form
$\delta(0)$ by $0$ and $\delta(\frak{n} + 1)$ by $\frak{n}$.  In later
stages, the set of terms is extended by function symbols introduced by
recursion, and the reduction of functionals there proceeds by
calculating the function for given numerical arguments according to
the recursive definition.  This will be discussed in the next section.
\end{enumerate}

In order to establish the consistency of the axiom system, Hilbert
suggests, we have to find a decidable (\emph{konkret feststellbar})
property of formulas so that every formula in a derivation which has
been transformed using the above steps has the property, and the
formula $0\neq 0$ lacks it. The property Hilbert proposes to use is
\emph{correctness}.  This is not to be understood as truth
in a model. The formulas still occurring in the derivation after the
transformation are all Boolean combinations of equations between numerals.
An equation between numerals $\frak{n} = \frak{m}$ is \emph{correct}
if $\frak{n}$ and $\frak{m}$ are syntactically equal, and the negation
of an equality is correct of $\frak{m}$ and $\frak{n}$ are not
syntactically equal.
\begin{quotation}
If we call a formula which does not contain variables or functionals
other than numerals an ``\emph{explicit [numerical] formula}'', then
we can express the result obtained thus: Every provable explicit
[numerical] formula is end formula of a proof all the formulas of
which are explicit formulas.

This would have to hold in particular of the formula $0 \neq 0$, if it
were provable. The required proof of consistency is thus completed if
we show that there can be no proof of the formula which consists of
only explicit formulas.

To see that this is impossible it suffices to find a concretely
determinable [\emph{konkret feststellbar}] property, which first of
all holds of all explicit formulas which result from an axiom by
substitution, which furthermore transfers from premises to end formula
in an inference, which however does not apply to the formula $0\neq
0$.\note{``\sg Nennen wir eine Formel, in der keine Variablen und
keine Funktionale ausser Zahlzeichen vorkommen, eine "`\emph{explizite
[numerische] Formel}"', so k"onnen wir das gefundene Ergebnis so
aussprechen: Jede beweisbare explizite [numerische] Formel ist
Endformel eines Beweises, dessen s"amtliche Formeln explizite Formeln
sind.

Dieses m"usste insbesondere von der Formel $0\neq 0$ gelten, wenn sie
beweisbar w"are. Der verlangte Nachweis der Widerspruchsfreiheit ist
daher erbracht, wenn wir zeigen, dass es keinen Beweis der Formel
geben kann, der aus lauter expliziten Formeln besteht.

Um diese Unm"oglichkeit einzusehen, gen"ugt es, eine konkret
feststellbare Eigenschaft zu finden, die erstens allen den expliziten
Formeln zukommt, welche durch Einsetzung aus einem Axiom entstehen,
die ferner bei einem Schluss sich von den Pr"amissen auf die Endformel
"ubertr"agt, die dagegen nicht auf die Formel $0\neq 0$ zutrifft.\eg''
\pcite{Hilbert:21}{part 2, 27--28}.}
\end{quotation}
Hilbert now defines the notion of a (conjunctive) normal form and
gives a procedure to transform a formula into such a normal form.  He
then provides the wanted property:
\begin{quotation}
With the help of the notion of a normal form we are now in a position
to exhipit a property which distinguishes the formula $0\neq 0$ from
the provable explicit formulas.

We divide the explicit formulas into ``\emph{correct}'' and
``\emph{incorrect}.'' The explicit atomic formulas are equations with
\emph{numerals} on either side [of the equality symbol].  We call such
an \emph{equation correct}, if the numerals on either side
\emph{coincide}, otherwise we call it \emph{incorrect}. We call an
\emph{inequality} with numerals on either side \emph{correct} if the
two numerals are \emph{different}, otherwise we call it
\emph{incorrect}.

In the normal form of an arbitrary explicit formula, each disjunct has
the form of an equation or an inequality with numerals on either side.

We now call a \emph{general explicit formula correct} if in the
corresponding normal form each disjunction which occurs as a conjunct
(or which constitutes the normal form) contains a correct equation or
a correct inequality. Otherwise we call the formula
\emph{incorrect}. [\dots]

According to this definition, the question of whether an explicit
formula is correct or incorrect is \emph{concretely decidable} in
every case.  Thus the ``tertium non datur'' holds here\dots\note{``\sg
Wir teilen die expliziten Formeln in "`\emph{richtige}"' und
"`\emph{falsche}"' ein. Die expliziten Primformeln sind Gleichungen,
auf deren beiden Seiten \emph{Zahlzeichen} stehen. Eine solche
\emph{Gleichung} nennen wir \emph{richtig}, wenn die beiderseits
stehenden Zahlzeichen \emph{"ubereinstimmen}; andernfalls nennen wir
sie \emph{falsch}.  Eine \emph{Ungleichung}, auf deren beiden Seiten
Zahlzeichen stehen, nennen wir \emph{richtig}, falls die beiden
Zahlzeichen \emph{verschieden} sind; sonst nen[n]en wir sie
\emph{falsch}.

In der Normalform einer beliebigen expliziten Formel haben alle
Disjunktionsglieder die Gestalt von Gleichungen oder Ungleichungen,
auf deren beiden Sieten Zahlzeichen stehen.

Wir nennen nun eine \emph{allgemeine explizite Formel richtig}, wenn
in der zugeh"origen Normalform jede als Konjunktionsglied auftretende
(bezw. die ganze Normalform ausmachende) Disjunktion eine richtige
Gleichung oder eine richtige Ungleichung als Glied
enth"alt. Andernfalls nennen wir die Formel \emph{falsch}. [\dots]

Nach der gegebenen Definition l"asst sich die Frage, ob eine
explizierte [sic] Formel richtig oder falsch ist, in jed[e]m Falle
\emph{konkret entscheiden}. Hier gilt also das "`tertium non
datur"' [\dots]\eg'' \pcite{Hilbert:21}{part~2, 33}.}
\end{quotation}
This use in the 1921--22 lectures of the conjunctive normal form of a
propositional formula to define correctness of Boolean combinations of
equalities between numerals goes back to the 1917--18 lecture
notes,\note{\pciteN{Hilbert:17}{149-150}. See also
\pciteN{Zach:99}{\S2.3}.} where transformation into conjunctive normal
form and testing whether each conjunct contains both $A$ and $\ol{A}$
was proposed as a test for propositional validity.  Similarly, here a
formula is \emph{correct} if each conjunct in its conjunctive normal
form contains a correct equation or a correct inequality.\note{A
sketch of the consistency proof is found in the Kneser
\emph{Mitschrift} to the 1921--22 lectures \cite{Hilbert:21a} in Heft
II, pp.~23--32 and in the official notes by Bernays
\pcite{Hilbert:21}{part 2, 19--38}.  The earlier Kneser
\emph{Mitschrift} leaves out step (1), and instead of eliminating
variables introduces the notion of \emph{einsetzungsrichtig}
(correctness under substitution, i.e., every substitution instance is
correct).  These problems were avoided in the official Bernays
typescript.  The Kneser notes did contain a discussion of recursive
definition and induction, which is not included in the official notes;
more about these in the next section.}  In the 1922--23 lectures, the
definition involving conjunctive normal forms is replaced by the usual
inductive definition of propositional truth and falsehood by truth
tables \pcite{Hilbert:22c}{21}.  Armed with the definition of correct
formula, Hilbert can prove that the derivation resulting from a proof
by transforming it according to (1)--(3) above contains only correct
formulas.  Since $0\neq 0$ is plainly not correct, there can be no
proof of $0\neq 0$ in the system consisting of axioms (1)--(10).  The
proof is a standard induction on the length of the derivation: the
formulas resulting from the axioms by elimination of variables and
reduction of functionals are all correct, and modus ponens preserves
correctness.\note{In the 1921--22 lectures, it is initially argued
that the result of applying transformations (1)--(3) results in a
\emph{proof} of the same end formula (if substitutions are added to
the initial formulas).  Specifically, it is suggested that the result
of applying elimination of variables and reduction of functionals to
the axioms results in formulas which are substitution instances of
axioms.  It was quickly realized that this is not the case. (When
Bernays presented the proof in the 1922--23 lectures on December~14,
1922, he comments that the result of the transformation need not be a
proof \pcite{Hilbert:22b}{21}.  The problem is the axiom of equality
\[
a = b \to (A(a) \to A(b)).
\]
Taking $A(c)$ to be $\delta(c) = c$, a substitution instance would be
\[
0+1+1 = 0 \to (\delta(0+1+1) = 0+1+1 \to \delta(0) = 0)
\]
This reduces to
\[
0+1+1 = 0 \to (0 + 1 = 0+1+1 \to 0 = 0)
\]
which is not a substitution instance of the equality axiom.  The
consistency proof itself is not affected by this, since the resulting
formula is still correct (in Hilbert's technical sense of the word).
The official notes to the 1921--22 lectures contain a 2-page
correction in Bernays's hand \pcite{Hilbert:21}{part~2, between pp.~26
and 27}.}

\subsection{Elementary number theory with recursion and induction rule}

The system of stage III consists of the basic system of the elementary
calculus of free variables and the successor function, extended by the
schema of defining functions by primitive recursion and the induction
rule.\note{The induction rule is not used in \cite{Ackermann:24},
since he deals with stage III only in passing and attempts a
consistency proof for all of analysis.  There, the induction rule is
superseded by an $\epsilon$-based induction axiom.  For a consistency
proof of stage III alone, an induction rule is needed, since an axiom
cannot be formulated without quantifiers (or $\epsilon$).  The
induction rule was introduced for stage III in the Kneser notes to the
1921--22 lectures (\citeA{Hilbert:21a}, Heft II, 32) and the 1922--23
lectures (\citeA{Hilbert:22b}, 26).  It is not discussed in the
official notes or the publications from the same period
(\citeA{Hilbert:22a}; \citeyear{Hilbert:23}} A primitive recursive
definition is a pair of axioms of the form
\begin{eqnarray*}
\phi(0, b_1, \ldots, b_n) & = & \frak{a}(b_1, \ldots, b_n)\\
\phi(a + 1, b_1, \ldots, b_n) & = & \frak{b}(a, \phi(a), b_1, \ldots, b_n)
\end{eqnarray*}
where $\frak{a}(b_1, \ldots, b_n)$ contains only the variables $b_1,$
\dots, $b_n$, and $\frak{b}(a, c, b_1, \ldots, b_n)$ contains only the
variables $a$, $c$, $b_1$, \ldots, $b_n$. Neither contains the
function symbol $\phi$ or any function symbols which have not yet been
defined.

The introduction of primitive recursive definitions and the induction
rule serves, first of all, the purpose of expressivity.  Surely any
decent axiom system for arithmetic must provide the means of
expressing basic number-theoretic states of affairs, and this includes
addition, subtraction, multiplication, division, greatest common
divisor, etc. The general schema of primitive recursion is already
mentioned in the Kneser notes for 1921--22 (\citeA{Hilbert:21a}, Heft
II, 29), and is discussed in some detail in the notes for the lectures
of the following year (\citeA{Hilbert:22c}, 26--30).

It may be interesting to note that in the 1922--23 lectures, there are
no axioms for addition or multiplication given before the general
schema for recursive definition.  This suggests a change in emphasis
during 1922, when Hilbert realized the importance of primitive
recursion as an arithmetical concept formation.  He later continued to
develop the notion, hoping to capture all number theoretic functions
using an extended notion of primitive recursion and to solve the
continuum problem with it.  This can be seen from the attempt at a
proof of the continuum hypothesis in \shortcite{Hilbert:26}, and
Ackermann's paper on ``Hilbert's construction of the reals''
\shortcite{Ackermann:28a}, which deals with hierarchies of recursive
functions.  The general outlook in this regard is also markedly
different from Skolem's \citeyear{Skolem:23}, which is usually
credited with the definition of primitive recursive
arithmetic.\note{The general tenor, outlook, and aims of Skolem's work
are sufficiently different from Hilbert's to suggest there was no
influence either way.  Skolem states in his concluding remarks that he
wrote the paper in 1919, after reading Russell and Whitehead's
\emph{Principia Mathematica}.  However, neither Hilbert nor Bernays's
papers contain an offprint or manuscript of Skolem's paper, nor
correspondence.  Skolem is not cited in any of Hilbert's, Bernays's,
or Ackermann's papers of the period, although the paper is referenced
in~\cite{HilbertBernays:34}.}

Hilbert would be remiss if he would not be including induction
in his arithmetical axiom systems.  As he already indicates in
the 1921--22 lectures, however, the induction principle cannot be 
formulated as an axiom without the help of quantifiers.  
\begin{quotation}
We are still completely missing the axiom of complete induction.  One
might think it would be
\[
\{Z(a) \to (A(a) \to A(a+1))\} \to \{A(1) \to (Z(b) \to A(b))\}
\] 
That is not it, for take $a=1$.  The hypothesis must hold for 
\emph{all $a$}.  We have, however, no means to bring the \emph{all}
into the hypothesis.  Our formalism does not yet suffice to 
write down the axiom of induction.

But as a schema we can: We extend our methods of proof by the following
schema.
\[
\infer{Z(a) \to \frak{K}(a)}{\frak{K}(1) & \frak{K}(a) \to \frak{K}(a+1)}
\]
Now it makes sense to ask whether this schema can lead to a
contradiction.\note{``\sg Uns fehlt noch ganz das
Ax[iom] der vollst["andigen] Induktion. Man k"onnte meinen, es w"are
\[
\{Z(a) \to (A(a) \to A(a+1))\} \to \{A(1) \to (Z(b) \to A(b))\}
\]
Das ist es nicht; denn man setze $a=1$.  Die Voraussetzung mu"s f"ur
\emph{alle $a$} gelten. Wir haben aber noch gar kein Mittel, das
\emph{Alle} in die Voraussetzung zu bringen. Unser Formalismus reich
noch nicht hin, das Ind.ax. aufzuschreiben.

Aber als Schema k"onnen wir es: Wir erweitern unsere Beweismethoden
durch das nebenstehende Schema
\[
\infer{Z(a) \to \frak{K}(a)}{\frak{K}(1) & \frak{K}(a) \to \frak{K}(a+1)}
\]
Jetzt ist es vern"unftig, zu fragen, ob dies Schema zum Wspruch
f"uhren kann.\eg'' \pcite{Hilbert:21a}{32}.  $Z$ is the predicate
expressing ``is a natural number,'' it disappears from later
formulations of the schema.}
\end{quotation}
The induction schema is thus necessary in the formulation of the
elementary calculus only because quantifiers are not yet available.
Subsequently, induction will be subsumed in the $\epsilon$-calculus.

The consistency proof for stage II is extended to cover also the
induction schema and primitive recursive definitions.  Both are only
sketched: Step (3), reduction of functionals, is extended to cover
terms containing primitive recursive functions by recursively
computing the value of the innermost term containing only numerals.
Both in the 1921--22 and the 1922--23 sets of notes by Kneser, roughly
a paragraph is devoted to these cases (the official sets of notes for
both lectures do not contain the respective passages).
\begin{quotation}
How de we proceed for recursions? Suppose a $\phi(\frak{z})$
occurs. Either [$\frak{z}$ is] 0, then we replace it by $\frak{a}$. Or
[it is of the form] $\phi(\frak{z+1})$: [replace it with]
$\frak{b}(\frak{z}, \phi(\frak{z}))$. Claim: These substitutions
eventually come to an end, if we replace innermost occurrences
first.\note{``\gq{Wie ist es bei Rekursionen? $\phi(\frak{z})$ komme
vor. Entweder 0, dann setzten wir $\frak{a}$ daf"ur. Oder
$\phi(\frak{z} + 1)$: $\frak{b}(\frak{z},
\phi\frak{z})$. Beh[auptung]: Das Einsetzen kommt zu einem Abschlu"s,
wenn wir zu innerst anfangen.}'' \pcite{Hilbert:22c}{29}}
\end{quotation}
The claim is not proved, and there is no argument that the process
terminates even for terms containing several different, nested
primitive recursively defined function symbols.

For the induction schema, Hilbert hints at how the consistency proof
must be extended.  Combining elimination of variables and reduction of
functionals we are to proceed upwards in the proof as before until we
arrive at an instance of the induction schema:
\[
\infer{Z(\frak{z}) \to \frak{K}'(\frak{z})}{
   \frak{K}(1) & \frak{K}(a) \to \frak{K}(a+1)}
\]
By copying the proof ending in the right premise, substituting
numerals $1$, \dots, $\frak{y}$ (where $\frak{z} = \frak{y}+1$) for
$a$ and applying the appropriate substitutions to the other variables
in $\frak{K}$ we obtain a proof of $Z(a) \to \frak{K}'(\frak{z})$
without the last application of the induction schema.

With the introduction of the $\epsilon$-calculus, the induction rule
is of only minor importance, and its consistency is never proved
in detail until \pciteN{HilbertBernays:34}{298--99}.

\subsection{The $\epsilon$-calculus and the axiomatization of mathematics}

In the spirit of the ``simultaneous development of logic and
mathematics,'' Hilbert takes the next step in the axiomatization of
arithmetic by employing a principle taken from Zermelo's
axiomatization of set theory: the axiom of choice.  Hilbert and
Bernays had dealt in detail with quantifiers in lectures in 1917--18
and 1920, but they do not directly play a significant role in the
axiom systems Hilbert develops for mathematics.  Rather, the first-
and higher-order calculi for which consistency proofs are proposed,
are based instead on choice functions.  The first presentation of
these ideas can be found in the 1921--22 lecture notes by Kneser (the
official notes do not contain these passages).  The motivation is that
in order to deal with analysis, one has to allow definitions of
functions which are not finitary.  These concept formations, necessary
for the development of mathematics free from intuitionist
restrictions, include definition of functions from undecidable
properties, by unbounded search, and choice.
\begin{quotation}
Not finitely (recursively) defined is, e.g., $\phi(a) = 0$ if there is
a $b$ so that $a^5 + ab^3 + 7$ is prime, and $=1$ otherwise.  But only
with these numbers and functions the real mathematical interest
begins, since the solvability in finitely many steps is not
foreseeable.  We have the conviction, that such questions, e.g., the
value of $\phi(a)$, are solvable, i.e., that $\phi(a)$ is also
finitely definable.  We cannot wait on this, however, we must allow
such definitions for otherwise we would restrict the free practice of
science.  We also need the concept of a function of
functions.\note{``\gq{Nicht endlich (durch Rek[ursion]) definiert ist
z.B. $\phi(a) = 0$ wenn es ein $b$ gibt, so da"s $a^5 + ab^3 + 7$
Primz[ahl] ist sonst $= 1$. Aber erst bei solchen Zahlen und
Funktionen beginnt das eigentliche math[ematische] Interesse, weil
dort die L"osbarkeit in endlich vielen Schritten nicht vorauszusehen
ist. Wir haben die "Uberzeugung, da"s solche Fragen wie nach dem Wert
$\phi(a)$ l"osbar, d.h. da"s $\phi(a)$ doch endlich definierbar
ist. Darauf k"onnen wir aber nicht warten: wir m"ussen solche
Definitionen zulassen, sonst w"urden wir den freien Betrieb der
Wissenschaft einschr"anken. Auch den Begriff der Funktionenfunktion
brauchen wir.}'' \pcite{Hilbert:21a}{Heft III, 1--2}.}
\end{quotation}
The concepts which Hilbert apparently takes to be fundamental for this
project are the principle of the excluded middle and the axiom of
choice, in the form of second-order functions $\tau$ and $\alpha$.
The axioms for these functions are
\begin{enumerate}
\item $\tau(f) = 0 \to (Z(a) \to f(a) = 1)$
\item $\tau(f) \neq 0 \to Z(\alpha(f))$
\item $\tau(f) \neq 0 \to f(\alpha(f)) \neq 1$
\item $\tau(f) \neq 0 \to \tau(f) = 1$
\end{enumerate}
The intended interpretation is: $\tau(f) = 0$ if $f$ is always $1$ and
$=1$ if one can choose an $\alpha(f)$ so that $f(\alpha(f)) \neq 1$.

The introduction of $\tau$ and $\alpha$ allows Hilbert to replace
universal and existential quantifiers, and also provides the 
basis for proofs of the axiom of induction and the least upper bound
principle.  Furthermore, Hilbert claims, the consistency of the
resulting system can be seen in the same way used to establish the
consistency of stage III (primitive recursive arithmetic). From a
proof of a numerical formula using $\tau$'s and $\alpha$'s, these
terms can be eliminated by finding numerical substitutions which turn
the resulting formulas into correct numerical formulas.

These proofs are sketched in the last part of the 1921--22 lecture
notes by Kneser \cite{Hilbert:21a}.\note{A full proof is given by
\citeN{Ackermann:24}.}  In particular, the consistency proof contains
the entire idea of the \emph{Hilbertsche Ansatz}, the
$\epsilon$-substitution method:
\begin{quotation}
First we show that we can eliminate all variables, since here also
only free variables occur.  We look for the innermost $\tau$ and
$\alpha$.  Below these there are only finitely defined [primitive
recursive] functions $\phi$, $\phi'$.  Some of these functions can be
substituted for $f$ in the axioms in the course of the proof.  1:
$\tau(\phi) = 0 \to (Z(\frak{a}) \to \phi(\frak{a})=1)$, where
$\frak{a}$ is a functional.  If this is \emph{not} used, we set all
$\alpha(\phi)$ and $\tau(\phi)$ equal to zero.  Otherwise we reduce
$\frak{a}$ and $\phi(\frak{a})$ and check whether $Z(\frak{a}) \to
\phi(\frak{a}) = 1$ is correct everywhere it occurs.  If it is
correct, we set $\tau[\phi] = 0$, $\alpha[\phi] = 0$.  If it is
incorrect, i.e., if $\frak{a} = \frak{z}$, $\phi(\frak{z}) \neq 1$, we
let $\tau(\phi) = 1$, $\alpha(\phi) = \frak{z}$.  After this
replacement, the proof remains a proof.  The formulas which take the
place of the axioms are correct.

(The idea is: if a proof is given, we can extract an argument from it
for which $\phi = 1$.)  In this way we eliminate the $\tau$ and
$\alpha$ and applications of [axioms] (1)--(4) and obtain a proof of
$1 \neq 1$ from I--V and correct formulas, i.e., from I--V,
\begin{eqnarray*}
\tau(f, b) = 0 & \to &\{Z(a) \to f(a, b) = 1\}\\
\tau(f, b) \neq 0 & \to &Z(f(\alpha, b))\\
\tau(f, b) \neq 0 & \to &f(\alpha(f, b), b) \neq 1\\
\tau(f, b) \neq 0 & \to &\tau(f, b) = 1\note{``\protect\sg
Als erstes zeigt man, da"s man alle Variablen fortschaffen kann, weil
auch hier nur freie Var[iable] vorkommen. Wir suchen die innersten
$\tau$ und $\alpha$. Unter diesen stehen nur endlich definierte
Funkt[ionen] $\phi$, $\phi'$\dots Unter diesen k"onnen einige im Laufe
des Beweises f"ur $f$ in die Ax[iome] eingesetzt sein.  1: $\tau(\phi)
= 0 \to (Z(\frak{a} \to \phi\frak{a} = 1)$ wo $\frak{a}$ ein Funktional
ist. Wenn dies \emph{nicht} benutzt wird, setze ich alle
$\alpha(\phi)$ und $\tau(\phi)$ gleich Null. Sonst reduziere ich
$\frak{a}$ und $\phi(\frak{a})$ und sehe, ob $Z(\frak{a} \to
\phi(\frak{a}) = 1$ in allen \dots wo sie vorkommt, richtig ist. Ist
die richtig, so setze ich $\tau = 0$ $\alpha = 0$. Ist sie falsch,
d.h. is $\frak{a} = \frak{z}$ $\phi(\frak{z} \neq 1$, so setzen wir
$\tau(\phi) = 1$, $\alpha(\phi) = \frak{z}$. Dabei bleibt der Beweis
Beweis. Die an Stelle der Axiome gesetzten Formeln sind richtig.

Der Gedanke ist: wenn ein Beweis vorliegt, so kann ich aus ihm ein
Argument finden f"ur das $\phi = 1$ ist). So beseitigt man
schrittweise die $\tau$ und $\alpha$ aund Anwendungen von 1 2 3 4 und
erh"alt einen Beweis von $1 \neq 1$ aus I--V und richtigen Formeln
d.h. aus I--V,\eg
\begin{eqnarray*}
\tau(f, b) = 0 & \to &\{Z(a) \to f(a, b) = 1\}\\
\tau(f, b) \neq 0 & \to &Z(f(\alpha, b))\\
\tau(f, b) \neq 0 & \to &f(\alpha(f, b), b) \neq 1\\
\tau(f, b) \neq 0 & \to &\tau(f, b) = 1\textrm{''}
\end{eqnarray*}
\cite{Hilbert:21a}, Heft III, 3--4.  The lecture is dated February
23, 1922.}
\end{eqnarray*}
\end{quotation}
Although not formulated as precisely as subsequent presentations, all
the ingredients of Hilbert's $\epsilon$-substitution method are here.
The only changes that are made en route to the final presentation of
Hilbert's sketch of the case for one $\epsilon$ and Ackermann's are
mostly notational.  In \cite{Hilbert:23}, a talk given in September
1922, the two functions $\tau$ and $\alpha$ are merged to one function
(also denoted $\tau$), which in addition provides the \emph{least}
witness for $\tau(f(a)) \neq 1$.  There the $\tau$ function is
also applied directly to formulas.  In fact, $\tau_a A(a)$ is the
primary notion, denoting the least witness $a$ for which $A(a)$ is
false; $\tau(f)$ is defined as $\tau_a(f(a) = 0)$.  Interestingly
enough, the sketch given there for the substitution method is for the
$\tau$-function for functions, not formulas, just as it was in the
1921--22 lectures.

The most elaborate discussion of the $\epsilon$-calculus can be found in 
Hilbert's and Bernays's course of 1922--23.  Here, again, the motivation
for the $\epsilon$-function is Zermelo's axiom of choice:
\begin{quotation}
What are we missing?
\begin{enumerate}
\item As far as logic is concerned: we have had the propositional
calculus extended by free variables, i.e., variables for which
arbitrary functionals may be substituted.  Operating with ``all'' and
``there is'' is still missing.
\item We have added the induction schema, but without consistency
proof and also on a provisional basis, with the intention of removing it.
\item So far only the arithmetical axioms which refer to whole numbers.
The above shortcomings prevent us from building up analysis (limit concept,
irrational number).
\end{enumerate}
These 3 points already give us a plan and goals for the following.

We turn to (1).  It is clear that a logic without ``all'---``there
is'' would be incomplete, I only recall how the application of these
concepts and of the so-called transfinite inferences has brought about
major problems.  We have not yet addressed the question of the
applicability of these concepts to infinite totalities.  Now we could
proceed as we did with the propositional calculus:  Formulate a few
simple [principles] as axioms, from which all others follow. Then
the consistency proof would have to be carried out---according to our
general program: with the attitude that a proof is a figure given to us.  
Significant obstacles to the consistency proof because of the bound 
variables.  The deeper investigation, however, shows that the real
core of the problem lies at a different point, to which one usually only
pays attention later, and which also has only been noticed in the
literature of late.\note{\sg``Was fehlt uns?
\begin{enumerate}
\item in logischer Hinsicht. Wir haben nur gehabt den Aussagenkalk"ul
mit der Erweiterung auf freie Variable d.h. solche f"ur die beliebige
Funktionale eingesetzt werden konnten.  Es fehlt das Operieren mit
"`alle"' und "`es gibt"'.
\item Wir haben das Induktionsschema hinzugef"ugt, ohne
W[iderspruchs]-f[reiheits] Beweis und auch nur provisorisch, also in
der Absicht, es wegzuschaffen.
\item Bisher nur die arithmethischen Axiome genau [?] die sich auf
ganze Zahlen beziehen.  Und die obigen M"angel verhindern uns ja
nat"urlich die Analysis aufzubauen (Grenzbegriff, Irrationalzahl).
\end{enumerate}
Diese 3 Punkte liefern schon Disposition und Ziele f"ur das Folgende.

Wir wenden uns zu 1.)  Es ist ja an sich klar, dass eine Logik ohne
"`alle"'---"`es gibt"' St"uckwerk w"are, ich erinnere wie gerade in
der Anwendung dieser Begriffe, und den sogennanten transfiniten
Schlussweisen die Hauptschwierigkeiten entstanden. Die Frage der
Anwendbarkeit dieser Begriffe auf $\infty$ Gesamtheiten haben wir noch
nicht behandelt. Nun k"onnten wir so verfahren, wie wir es beim
Aussagen-Kalk"ul gemacht haben: einige, m"oglichst einfache [???] als
Axiome zu formalisieren, aus denen sich [sic] dann alle "ubrigen
folgen.  Dann m"usste der W-f Beweis gef"uhrt werden---unserem
allgem[einen] Programm gem"ass: mit unserer Einstellung, dass Beweis
eine vorliegende Figur ist.  F"ur den W-f Beweis grosse
Schwierigkeiten wegen der gebundenen Variabeln.  Die tiefere
Untersuchung zeigt aber, dass der eigentliche Kern der Schwierigkeit
an einer anderen Stelle liegt, auf die man gew"ohnlich erst sp"ater
Acht giebt und die auch in der Litteratur erst sp"ater wahrgenommen
worden ist.''\eg \pcite{Hilbert:22b}{\emph{Erg"anzung}, sheets 3--4}.}
\end{quotation}
At the corresponding place in the Kneser \emph{Mitschrift}, Hilbert
continues:
\begin{quotation}
[This core lies] in Zermelo's \emph{axiom of choice}. [\dots] The
objections [of Brouwer and Weyl] are directed against the choice
principle. But they should likewise be directed against ``all'' and
``there is'', which are based on the same basic idea.

We want to extend the axiom of choice.  To each proposition
with a variable $A(a)$ we assign an object for which the
proposition holds only if is holds in general.  So, a counterexample,
if one exists.

$\epsilon(A)$, an individual logical function. [\dots] $\epsilon$
satisfies the \emph{transfinite axiom}:
\begin{enumerate}
\item[(16)] $A(\epsilon A) \to Aa$
\end{enumerate}
e.g., $Aa$ means: $a$ is corrupt. $\epsilon A$ is Aristides.\note{\sg 
``[Dieser Kern liegt] beim \emph{Auswahlaxiom} von Zermelo. [\dots] Die
Einw"ande richten sich gegen das Auswahlprinzip. Sie m"u"sten sich
aber ebenso gegen "`alle"' und "`es gibt"' richten, wobei derselbe
Grundgedanke zugrunde liegt.

Wir wollen das Auswahlaxiom erweitern.  Jeder Aussage mit einer
Variablen $A(a)$ ordnen wir ein Ding zu, f"ur das die Aussage nur dann
gilt, wenn sie allgemein gilt.  Also ein Gegenbeispiel, wenn es
existiert.

$\epsilon(A)$, eine individuelle logische Funktion. [\dots] $\epsilon$
gen"uge dem \emph{transfiniten Axiom}:\medskip

\noindent (16)\qquad $A(\epsilon A) \to Aa$\medskip

\noindent z.B. $Aa$ hei"se: $a$ ist bestechlich. $\epsilon A$ ist
Aristides.''\eg \pcite{Hilbert:22c}{30--31}.  The lecture is dated
February~1, 1922, given by Hilbert.  The corresponding part of
Hilbert's notes for that lecture in \pcite{Hilbert:22b}{Erg"anzung,
sheet~4}) contains page references to \pcite{Hilbert:23}{152 and 156,
paras. 4--6 and 17--19 of the English translation}, and indicates the
changes made for the lecture, specifically, to replace $\tau$ by
$\epsilon$.}
\end{quotation}

Hilbert goes on to show how quantifiers can be replaced by
$\epsilon$-terms.  The corresponding definitional axioms are already
included in \citeN{Hilbert:23}, i.e., $A(\epsilon A) \equiv (a)A(a)$
and $A(\epsilon \ol{A}) \equiv (\exists a)A(a)$.  Next, Hilbert outlines
a derivation of the induction axioms using the $\epsilon$-axioms.  For
this, it is necessary to require that the choice function takes the
minimal value, which is expressed by the additional axiom
\[
\epsilon A \neq 0 \to A(\delta(\epsilon A)).
\]
With this addition, Hilbert combined the $\kappa$ function of
\cite{Hilbert:22a} and the $\mu$ function of \shortcite{Hilbert:23}
with the $\epsilon$ function.  Both $\kappa$ (``k'' for
\emph{Kleinstes}, least) and $\mu$ had been introduced there as
functions of functions giving the least value for which the function
differs from $0$. In \pciteN{Hilbert:23}{161--162}, Hilbert indicates
that the axiom of induction can be derived using the $\mu$ function,
and credits this to \citeN{Dedekind:88}.

The third issue Hilbert addresses is that of dealing with real
numbers, and extending the calculus to analysis.  A first step can be
carried out at stage IV by considering a real number as a function
defining an infinite binary expansion.  A sequence of reals can then
be given by a function with two arguments.  Already in
\citeN{Hilbert:23} we find a sketch of the proof of the least upper
bound principle for such a sequence of reals, using the $\pi$
function:
\[
\pi A(a) = \cases{0 & if $(a)A(a)$ \cr 1 & otherwise}
\]
The general case of sets of reals needs function variables and
second-order $\epsilon$ and $\pi$. These are briefly introduced as
$\epsilon_f A$ with the axioms
\begin{eqnarray*}
A\epsilon_f A &\to & Af \\
(f)Af & \to & \pi_f Af = 0\\
\ol{(f)}Af & \to & \pi_f Af = 1
\end{eqnarray*}
The last two lectures transcribed in \cite{Hilbert:22c} are devoted to
a sketch of the $\epsilon$ substitution method.  The proof is adapted
from \citeN{Hilbert:23}, replacing $\epsilon f$ with $\epsilon A$,
also deals with $\pi$, and covers the induction axiom in its
form for the $\epsilon$-calculus.\note{See section~\ref{sec-esub} on
the $\epsilon$-substitution method.}  During the last lecture, Bernays
also extends the proof to second-order $\epsilon$'s.
\begin{quotation}
If we have a \emph{function variable}:
\[
A \epsilon_f Af \to Af
\]
[\dots] Suppose $\epsilon$ only occurs with $\frak{A}$ (e.g., $f0=0$,
$ff0 = 0$).  How will we eliminate the function variables?  We simply
replace $fc$ by $c$.  This does not apply to \emph{bound} variables.
For those we take some fixed function, e.g., $\delta$ and carry out
the reduction with it.  Then we are left with, e.g., $\frak{A}\delta
\to \frak{A}\phi$. This, when reduced, is either correct or incorrect.  In
the latter case, $\frak{A}\phi$ is incorrect.  Then we substitute $\phi$
everywhere for $\epsilon_f \frak{A}f$.  Then we end up with
$\frak{A}\phi \to \frak{A}\psi$.  That is certainly correct, since
$\frak{A}\phi$ is incorrect.\note{``\sg Wenn wir eine
\emph{Funktionsvariable} haben:
\[
A \epsilon_f Af \to Af
\]
($\pi$ f"allt fort)? $\epsilon$ komme \emph{nur} mit $\frak{A}$ vor
(z.B. $f0 = 0$, $ff0 = 0$).  Wie werden wir die Funktionsvariablen
ausschalten?  Statt $fc$ setzen wir einfach $c$. Auf die
\emph{gebundenen} trifft das nicht zu.  F"ur diese nehmen wir
probeweise eine bestimmte Funktion z.B. $\delta a$ und f"uhren damit
die Reduktion durch.  Dann steht z.B. $\frak{A}\delta \to
\frak{A}\phi$. Diese reduziert ist r[ichtig] oder f[alsch]. Im letzten
Falle is $\frak{A}\phi$ falsch. Dann setzen wir "uberall $\phi$ f"ur
$\epsilon_f\frak{A}f$.  Dann steht $\frak{A}\phi \to
\frak{A}\psi$. Das ist sicher \emph{r[ichtig]} da $\frak{A}\phi$
f[alsch] ist.\eg'' \pcite{Hilbert:22c}{38--39}.}
\end{quotation}

The last development regarding the $\epsilon$-calculus before
Ackermann's dissertation is the switch to the dual notation.  Both
\cite{Hilbert:23} and \cite{Hilbert:22c} use $\epsilon A$ as denoting
a counterexample for $A$, whereas at least from Ackermann's
dissertation onwards, $\epsilon A$ denotes a witness.
Correspondingly, Ackermann uses the dual axiom $A(a) \to A(\epsilon_a
A(a))$.  Although it is relatively clear that the supplement to
\cite{Hilbert:22c}---24 sheets in Hilbert's hand---are Hilbert's notes
based on which he and partly Bernays presented the 1922--23 lectures,
parts of it seem to have been altered or written after the conclusion
of the course. Sheets 12--14 contain a sketch of the proof of the
axiom of induction from the standard, dual $\epsilon$ axioms; the same
proof for the original axioms can be found on sheets 8--11.

This concludes the development of mathematical systems using the
$\epsilon$-calculus and consistency proofs for them presented by
Hilbert himself.  We now turn to the more advanced and detailed
treatment in Wilhelm Ackermann's \shortcite{Ackermann:24}
dissertation.

\section{Ackermann's Dissertation}

Wilhelm Ackermann was born in 1896 in Westphalia.  He studied
mathematics, physics, and philosophy in G\"ottingen between 1914 and
1924, serving in the army in World War I from 1915--1919.  He
completed his studies in 1924 with a dissertation, written under
Hilbert, entitled ``\gq{Begr"undung des `tertium non datur' mittels
der Hilbertschen Theorie der Widerspruchsfreiheit}''
(\citeA{Ackermann:24a}; \citeyear{Ackermann:24}), the first major
contribution to proof theory and Hilbert's Program.  In 1927 he
decided for a career as a high school teacher rather than a career in
academia, but remained scientifically active.  His major contributions
to logic include the function which carries his name---an example of a
recursive but not primitive recursive function \cite{Ackermann:28a},
the consistency proof for arithmetic using the $\epsilon$-substitution
method \cite{Ackermann:40}, and his work on the decision problem
(\citeA{Ackermann:28}; \citeyear{Ackermann:54}).  He served as
co-author, with Hilbert, of the influential logic textbook
\emph{\gq{Grundz"uge der theoretischen Logik}} \cite{Hilbert:28b}.  He
died in 1962.\note{For a more detailed survey of Ackermann's
scientific contributions, see \cite{Hermes:67}.  A very informative
discussion of Ackermann's scientific correspondence can be found in
\cite{Ackermann:83}.}

Ackermann's 1924 dissertation is of particular interest since it is
the first non-trivial example of what Hilbert considered to be a
finitistic consistency proof.  Von Neumann's paper of 1927 does not
entirely fit into the tradition of the Hilbert school, and we have no
evidence of the extent of Hilbert's involvement in its writing.  Later
consistency proofs, in particular those by Gentzen and Kalm\'ar, were
written after G\"odel's incompleteness results were already well-known
and their implications understood by proof theorists.  Ackermann's
work, on the other hand, arose entirely out of Hilbert's research
project, and there is ample evidence that Hilbert was aware of the
range and details of the proof.  Hilbert was Ackermann's dissertation
advisor, approved the thesis, was editor of {\em Mathematische
Annalen,} where the thesis was published, and corresponded with
Ackermann on corrections and extensions of the result.  Ackermann was
also in close contact with Paul Bernays, Hilbert's assistant and close
collaborator in foundational matters.  Ackermann spent the first half
of 1925 in Cambridge, supported by a fellowship from the International
Education Board (founded by John D. Rockefeller, Jr., in 1923).  In
his letter of recommendation for Ackermann, Hilbert writes:
\begin{quotation}
In his thesis ``Foundation of the `tertium non datur' using Hilbert's
theory of consistency,'' Ackermann has shown in the most general case that
the use of the words ``all'' and ``there is,'' of the ``tertium non
datur,'' is free from contradiction. The proof uses exclusively
primitive and finite inference methods.  Everything is demonstratded,
as it were, directly on the mathematical formalism.

Ackermann has here surmounted considerable mathematical difficulties
and solved a problem which is of first importance to the modern
efforts directed at providing a new foundation for
mathematics.\note{``\sg In seiner Arbeit ``Begr"undung des "`Tertium
non datur"' mittels der Hilbertschen Theorie der Widerspruchsfreiheit
hat Ackermann im allgemeinsten Falle gezeigt, dass der Gebrauch der
Worte "`alle"' und "`es gibt"', des "`Tertium non datur"'
widerspruchsfrei ist. Der Beweis erfolgt unter ausschliesslicher
Benutzung primitiver und endlicher Schlussweisen. Es wird alles an dem
mathematischen Formalismus sozusagen direkt demonstriert.

Ackermann hat damit unter Ueberwindung erheblicher mathematischer
Schwierigkeiten ein Problem gel"ost, das bei den modernene auf eine
Neubegr"undung der Mathematik gerichteten Bestrebungen an erster
Stelle steht.\eg'' Hilbert-Nachla"s, Nie\-der\-s"achsische Staats- und
Universit"atsbibliothek G"ottingen, Cod. Ms. Hilbert 458, sheet 6, no
date.  The three-page letter was evidently written in response to a
request by the President of the International Education Board, dated
May 1, 1924.}
\end{quotation}
Further testimony of Hilbert's high esteem for Ackermann can be found
in the draft of a letter to Russell asking for a letter of support to
the International Education Board, where he writes that ``Ackermann
has taken my classes on foundations of mathematics in recent semesters
and is currently one of the best masters of the theory which I have
developed here.''\note{``\gq{Ich bemerke nur, dass Ackermann meine
Vorlesungen "uber die Grundlagen der Math\-[ematik] in den letzten
Semestern geh"ort hat und augenblicklich einer der besten Herren der
Theorie ist, die ich hier entwickelt habe.}'' \emph{ibid.}, sheet 2.
The draft is dated March~19, 1924, and does not mention Russell by
name.  \citeN{Sieg:99}, however, quotes a letter from Russell's wife
to Hilbert dated May~20, 1924, which responds to an inquiry by
Hilbert concerning Ackermann's stay in Cambridge.  Later in the
letter, Hilbert expresses his regret that the addressee still has not
been able to visit G"ottingen.  Sieg documents Hilbert's effort in the
preceding years to effect a meeting in G"ottingen; it is therefore
quite likely that the addressee was Russell.}

Ackermann's work provides insight into two important issues relating
to Hilbert's program as it concerns finitistic consistency proofs.
First, it provides historical insight into the aims and development of
Hilbert's Program: The first part of the program called for an
axiomatization of mathematics.  These axiomatizations were then the
objects of metamathematical investigations: the aim was to find
finitistic consistency proofs for them.  Which areas of mathematics
were supposed to be covered by the consistency proofs, how were they
axiomatized, what is the strength of the systems so axiomatized?  We
have already seen what Hilbert's roadmap for the project of
axiomatization was.  Ackermann's dissertation provides the earliest
example of a formal system stronger than
elementary arithmetic.  The second aim, the metamathematical
investigation of the formal systems obtained, also poses historical
questions: When did Ackermann, and other collaborators of Hilbert (in
particular, Bernays and von Neumann) achieve the results they sought?
Was Ackermann's proof correct, and if not, what parts of it can be
made to work?

The other information we can extract from an analysis of Ackermann's
work is what methods were used or presupposed in the consistency
proofs that were given, and thus, what methods were sanctioned by
Hilbert himself as falling under the finitist standpoint.  Such an
analysis of the methods used are of a deeper, conceptual interest.
There is a fundamental division between Hilbert's philosophical
remarks on finitism on the one hand, and the professed goals of the
program on the other.  In these comments, rather little is said about
the concept formations and proof methods that a finitist, according to
Hilbert, is permitted to use.  In fact, most of Hilbert's remarks deal
with the objects of finitism, and leave the finitistically admissible
forms of definition and proof to the side.  These, however, are the
questions at issue in contemporary conceptual analyses of
finitism. Hilbert's relative silence on the matter is responsible for
the widespread---and largely correct---opinion that Hilbert was too
vague on the question of what constitutes finitism to unequivocally
define the notion, and therefore later commentators have had enough
leeway to disagree widely on the strength of the finitist standpoint
while still claiming to have explicated Hilbert's own concept.

\subsection{Second-order primitive recursive arithmetic}\label{sec-2pra}

In \citeN{Ackermann:24}, the system of stage III is extended by
second-order variables for functions.  The schema of recursive
definition is then extended to include terms containing such variables.
In the following outline, I shall follow Ackermann and adopt the notation of
subscripting function symbols and terms by variables to indicate that
these variables do not occur freely but rather as placeholders for
functions.  For instance, $\frak{a}_a(f(a))$ indicates that the term
$\frak{a}$ does not contain the variable $a$ free,  but rather that 
the function $f(a)$ appears as a functional argument, i.e., 
that the term is of the form $\frak{a}(\lambda a.f(a))$.
The schema of second-order primitive recursion is the following:
\begin{eqnarray*}
\phi_{\vec b_i}(0, \vec{f}(\vec{b}_i), \vec{c}) & = & 
\frak{a}_{\vec{b}_i}(\vec{f}(\vec{b}_i), \vec{c})\\
\phi_{\vec{b}_i}(a+1, \vec f(\vec{b}_i), \vec c) & = & 
\frak{b}_{\vec{b}_i}(a, \phi_{\vec{d}_i}(a, \vec f(\vec{d}_i), \vec c), \vec f(\vec{b}_i))
\end{eqnarray*}
To clarify the subscript notation, compare this with the schema of 
second-order primitive recursion using $\lambda$-abstraction notation:
\begin{eqnarray*}
\phi(0, \lambda \vec{b}_i.\vec{f}(\vec{b}_i), \vec{c}) & = & 
\frak{a}(\lambda \vec{b}_i.\vec{f}(\vec{b}_i), \vec{c})\\
\phi(a+1, \lambda \vec{b}_i.\vec f(\vec{b}_i), \vec c) & = & 
\frak{b}(a, \phi(a, \lambda \vec{d}_i.\vec f(\vec{d}_i), \vec c), \lambda \vec{b}_i.\vec f(\vec{b}_i))
\end{eqnarray*}
Using this schema, it is possible to define the Ackermann function.
This was already pointed out in \citeN{Hilbert:26}, although it was
not until \citeN{Ackermann:28a} that it was shown that the function so
defined cannot be defined by primitive recursion without function
variables.  \citeN{Ackermann:28a} defines the function as follows.
First it is observed that the iteration function 
\[
\rho_c(a, f(c), b) = \underbrace{f(\dots f(f}_{\textrm{$a$ $f$'s}}(b))\dots)
\]
can be defined by second-order primitive recursion:
\begin{eqnarray*}
\rho_c(0, f(c), b) & = & b \\
\rho_c(a +1, f(c), b) & = & f(\rho_c(a, f(c), b))
\end{eqnarray*}
Furthermore, we have two auxiliary functions 
\[
\iota(a, b) = \cases{1 & if $a = b$\cr 0 & if $a \neq b$} \quad \textrm{and}\quad
\lambda(a, b) = \cases{0 & if $a \neq b$\cr 1 & if $a = b$}
\]
which are primitive recursive, as well as addition and multiplication.
The term $\frak{a}(a, b)$ is short for $\iota(a, 1)\cdot \iota(a, 0) \cdot b + \lambda(a, 1)$; we then have
\[
\frak{a}(a, b) = \cases{0 & if $a =0$\cr 1 & if $a =1$ \cr b & otherwise}
\]
The Ackermann function is defined by
\begin{eqnarray*}
\phi(0, b, c) & = & b + c\\
\phi(a+1, b, c) & = & \rho_d(c, \phi(a, b, d), \frak{a}(a, b)).
\end{eqnarray*}
In more suggestive terms,
\begin{eqnarray*}
\phi(0, b, c) & = & b + c\\
\phi(1, b, c) & = & b \cdot c\\
\phi(a+1, b, c) &= &\underbrace{\phi(a, b, \phi(a, b, \ldots 
\phi(a, b, b)\cdots))}_{c \textrm{\ times}} 
\end{eqnarray*}
The system of second-order primitive recursive arithmetic 2PRA$^-$
used in \citeN{Ackermann:24} consists of axioms (1)--(15) of
\pciteN{Hilbert:22b}{see Section~\ref{sec-early-conprfs}}, axiom (16)
was replaced by
\begin{enumerate}
\item[16.] $a \neq 0 \to a = \delta(a) + 1$,
\end{enumerate}
plus defining equations for both first- and second-order primitive
recursive functions.  There is no induction rule (which is usually
included in systems of primitive recursive arithmetic), although the
consistency proof given by Ackermann can easily be extended to cover
it.


\subsection{The consistency proof for primitive recursive arithmetic}

Allowing primitive recursion axioms for functions which contain
function variables is a natural extension of the basic calculus of
stages III and IIII, and is necessary in order to be able to introduce
sufficiently complex functions.  Hilbert seems to have thought that by
extending primitive recursion in this way, or at least by building an
infinite hierarchy of levels of primitive recursions using variables
of higher types, he could account for \emph{all} the number theoretic
functions, and hence for all real numbers (represented as decimal
expansions).  In the spirit of the stage-by-stage development of
systems of mathematics and consistency proofs, it is of course
necessary to show the consistency of the system of stage IIII, which
is the system presented by Ackermann.  As before, it makes perfect
sense to first establish the consistency for the fragment of stage
IIII not containing the transfinite $\epsilon$ and $\pi$ functions.
In Section 4 of his dissertation, Ackermann undertakes precisely this
aim.

The proof is a direct extension of the consistency proof of stage III,
the elementary calculus of free variables with basic arithmetical
axioms and primitive recursive definitions, i.e., PRA.  This proof had
already been presented in Hilbert's lectures in 1921--22 and 1922--23.
The idea here is the same: put a given, purported proof of $0 \neq 0$
into tree form, eliminate variables, and reduce functionals.  The
remaining figure consists entirely of correct formulas, where correctness
of a formula is a syntactically defined and easily decidable property.
The only complication for the case where function variables are also
admitted is the reduction of functionals.  It must be shown that every
functional, i.e., every term of the language, can be reduced to a
numeral on the basis of the defining recursion equations.  For the
original case, this could be done by a relatively simple inductive
proof.  For the case of 2PRA$^-$, it is not so obvious.

Ackermann locates the difficulty in the following complication.
Suppose you have a functional $\phi_b(2, \frak{b}(b))$.  Here,
$\frak{b}(b)$ is a term which denotes a function, and so there is no
way to replace the variable $b$ with a numeral before evaluating the
entire term.  In effect, the variable $b$ is bound (in modern
notation, the term might be more suggestively written $\phi(2, \lambda
b.\frak{b}(b))$.)  In order to reduce this term, we apply the
recursion equations for $\phi$ and might end up with a
term like
\[
\frak{b}(1) + \frak{b}(2) + \frak{b}(0)\cdot\frak{b}(1) + 
\frak{b}(1)\cdot\frak{b}(2).
\]
The remaining $\frak{b}$'s might in turn contain $\phi$.  By contrast,
reducing a term $\psi(\frak{z})$ where $\psi$ is defined by
first-order primitive recursion results in a term which does not
contain $\psi$. 

To show that the reduction indeed comes to an end if innermost
subterms are reduced first, Ackermann proposes to assign indices to
terms and show that each reduction reduces this index.  The indices
are, essentially, ordinal notations $< \omega^{\omega^\omega}$.  Since
this is probably the first proof using ordinal notations, it may be of
some interest to repeat and analyze it in some detail here.  In my
presentation, I stay close to Ackermann's argument and only change the
notation for ranks, indices, and orders: Where Ackermann uses sequences
of natural numbers, I will use the more perspicuous ordinal notations.

Suppose the primitive recursive functions are arranged in a linear
order according to the order of definition.  We write $\phi < \psi$ if
$\phi$ occurs before $\psi$ in the order of definition, i.e., $\psi$
cannot be used in the defining equations for $\phi$.  Suppose further
that we are given a specific term $t$.  The notion of
\emph{subordination} is defined as follows: an occurrence of a
function symbol $\xi$ in $t$ is subordinate to an occurrence of
$\phi$, if $\phi$ is the outermost function symbol of a subterm $s$,
the occurrence of $\xi$ is in $s$, and the subterm of $s$ with
outermost function symbol $\xi$ contains a bound variable $b$ in whose
scope the occurrence of $\phi$ is (this includes the case where $b$
happens to be bound by $\phi$ itself).\note{\label{change}Ackermann
only requires that $b$ be bound by the occurrence of $\phi$, but this
is not enough for his proof.}  In other words, $\frak{t}$ is of the form
\[
\frak{t}'(\ldots \phi_b(\ldots \xi(\ldots b
\ldots)\ldots)\ldots)
\]
The \emph{rank} $rk(\frak{t}, \phi)$ of an \emph{occurrence} of a function
symbol $\phi$ with respect to $\frak{t}$ is defined as follows: If there is
no occurrence of $\psi > \phi$ which is subordinate to $\phi$ in $\frak{t}$,
then $rk(\frak{t}, \phi) = 1$.  Otherwise,
\[ 
rk(\frak{t}, \phi) = \max\{rk(\frak{t}, \psi) : \mbox{$\psi > \phi$ is subordinate
to $\phi$}\} + 1.
\]
The rank $r(\frak{t}, \phi)$ of a term $\frak{t}$ with respect to a
function symbol $\phi$ is the maximum of the ranks of occurrences of
$\phi$ or $\psi > \phi$ in $t$.  (If neither $\phi$ nor $\psi > \phi$
occur in $\frak{t}$, that means $r(\frak{t}, \phi) = 0$.

We assign to a subterm $\frak{s}$ of $\frak{t}$ a sequence of ranks of
$\psi_n$, \ldots, $\psi_0$ with respect to $\frak{s}$, where $\psi_0 <
\cdots < \psi_n$ are all function symbols occurring in $\frak{t}$.
This is the \emph{order} of~$\frak{s}$:
\[
o(\frak{s}) = \langle r(\frak{s}, \psi_n), \ldots, r(\frak{s}, \psi_0)\rangle
\]
We may think of this order as an ordinal $< \omega^\omega$, specifically,
$o(\frak{s})$ corresponds to
\[
\omega^n \cdot r(\frak{s}, \psi_n)\cdot + \cdots + \omega\cdot r(\frak{s},
\psi_1)\cdot + r(\frak{s}, \psi_0)
\]
Now consider the set of all distinct subterms of $t$ of a given order
$o$ which are not numerals. The \emph{degree} $d(\frak{t}, o)$ of $o$
in $\frak{t}$ is the cardinality of that set.  The \emph{index} $j(\frak{t})$
of $\frak{t}$ is the sequence of degrees ordered in the same way as the
orders, i.e., it corresponds to an ordinal of the form
\[
\sum_o \omega^o \cdot d(\frak{t}, o)
\]
where the sum ranges over all orders $o$.  Obviously, this is an
ordinal $< \omega^{\omega^\omega}$.

Ackermann does not use ordinals to define these indices; he stresses
that he is only dealing with finite sequences of numbers, on which an
elaborate order is imposed.  Rather than appeal to the
well-orderedness of $\omega^{\omega^\omega}$, he gives a more direct
argument that by repeatedly proceeding to indices which are smaller in
the imposed order one eventually has to reach the index which consists
of all~$0$.

Suppose a term $\frak{t}$ not containing $\epsilon$ or $\pi$ is given.
Let $\frak{s}$ be an innermost constant subterm which is not a
numeral, we may assume it is of the form $\phi_{\bar b}(\frak{z}_1,
\ldots, \frak{z}_n, \frak{u}_1, \ldots, \frak{u}_m)$ where
$\frak{u}_i$ is a term with at least one variable bound by $\phi$ and
which doesn't contain a constant subterm.  We have two cases:

(1) $\frak{s}$ does not contain bound variables, i.e., $m=0$.  The
order of $\frak{s}$ is $\omega^k$ (where $\phi = \psi_k$).  Evaluating
the term $\frak{s}$ by recursion results in a term $\frak{s}'$ in
which only function symbols of lower index occur. Hence, the highest
exponent in the order of $\frak{s}'$ is less than $i$, and so $o(s')
<o(s)$.  Furthermore, since no variable which is bound in $t$ can
occur in $\frak{s}'$ (since no such variable occurs in $\frak{s}$),
replacing $\frak{s}$ by $\frak{s}'$ in $\frak{t}$ does not result in
new occurrences of function symbols which are subordinate to any
other.  Thus the number of subterms in the term $\frak{t}'$ which
results from such a replacement with orders $> o(\frak{s})$ remains
the same, while the number of subterms of order $o(\frak{s})$ is
reduced by~1.  Hence, $j(\frak{t}') < j(\frak{t})$.

(2) $\frak{s}$ does contain bound variables.  For simplicity, assume that
there is one numeral argument and one functional argument, i.e.,
$\frak{s}$ is of the form $\phi_b(\frak{z}, \frak{c}(b))$.  In this
case, all function symbols occurring in $\frak{c}(b)$ are subordinate
to $\phi$, or otherwise $\frak{c}(b)$ would contain a constant
subterm.\note{According to Ackermann's definition of subordination,
this would not be true.  A subterm of $\frak{c}(b)$ might contain a
bound variable and thus not be a constant subterm, but the variable
could be bound by a function symbol in $\frak{t}$ other than the
occurrence of $\phi$ under consideration.  See note~\ref{change}.}
Thus, the rank of $\frak{c}(b)$ in $\frak{t}$ with respect to $\psi_i$
is less that the rank of $\frak{s}$ with respect to $\psi_i$.

We reduce the subterm $\frak{s}$ to a subterm $\frak{s'}$ by applying
the recursion.  $\frak{s'}$ does not contain the function symbol
$\phi$.  We want to show that replacing $\frak{s}$ by $\frak{s}'$ in
$\frak{t}$ lowers the index of $\frak{t}$.

First, note that when substituting a term $\frak{a}$ for $b$ in
$\frak{c}(b)$, the order of the resulting $\frak{c}(\frak{a})$ with
respect to $\phi$ is the maximum of the orders of $\frak{c}(b)$ and
$\frak{a}$, since none of the occurrences of function symbols in
$\frak{a}$ contain bound variables whose scope begins outside of
$\frak{a}$, and so none of these variables are subordinate to any
function symbols in $\frak{c}(b)$.  

Now we prove the claim by induction on $\frak{z}$.  Suppose the defining 
equation for $\phi$ is
\begin{eqnarray*}
\phi_b(0, f(b)) & = & \frak{a}_b(f(b)) \\
\phi_b(a+1, f(b)) & = & \frak{b}_b(\phi_c(a, f(c)), a, f(b)).
\end{eqnarray*}
If $\frak{z} = 0$, then $\frak{s}' = \frak{a}_b(\frak{c}(b))$.  At a
place where $f(b)$ is an argument to a function, $f(b)$ is replaced by
$\frak{c}(d)$, and $d$ is not in the scope of any $\phi$ (since
$\frak{a}$ doesn't contain $\phi$).  For instance, $\frak{a}_b(f(b)) =
2 + \psi_d(3, f(d))$. Such a replacement cannot raise the $\phi$-rank
of $\frak{s'}$ above that of $\frak{c}(b)$.  The term $\frak{c}$ might
also be used in places where it is not a functional argument, e.g., if
$\frak{a}_b(f(b)) = f(\psi(f(2)))$.  By a simple induction on the
nesting of $f$'s in $\frak{a}_b(f(b))$ it can be seen that the
$\phi$-rank of $\frak{s'}$ is the same as that of $\frak{c}(b)$: For
$\frak{c}(\frak{d})$ where $\frak{d}$ does not contain $\frak{c}$,
the $\phi$-rank of $\frak{c}(\frak{d})$ equals that of
$\frak{c}(b)$ by the note above and the fact that $\frak{d}$ does not
contain $\phi$.  If $\frak{d}$ does contain a nested $\frak{c}$, then
by induction hypothesis and the first case, its $\phi$-rank is the
same as that of $\frak{c}(b)$.  By the note, again, the entire subterm
has the same $\phi$-rank as $\frak{c}(b)$.

The case of $\phi_b(\frak{z}+1, \frak{c}(b))$ is similar.  Here, the
first replacement is \[\frak{b}_b(\phi_c(\frak{z}, \frak{c}(c)),
\frak{z}, \frak{c}(b)).\]  Further recursion replaces $\phi_c(\frak{z},
\frak{c}(c))$ by another term which, by induction hypothesis, has
$\phi$-rank less than or equal to that of $\frak{c}(b)$.  The same
considerations as in the base case show that the entire term also has
a $\phi$-rank no larger than $\frak{c}(b)$.

We have thus shown that eliminating the function symbol $\phi$ by
recursion from an innermost constant term reduces the $\phi$-rank of
the term at least by one and does not increase the $\psi_j$-ranks of
any subterms for any $j > i$.

In terms of ordinals, this shows that at least one subterm of order $o$
was reduced to a subterm of order $o' < o$, all newly introduced  subterms have
order $< o$, and the order of no old subterm increased.  Thus, the index
of the entire term was reduced.  The factor $\omega^o \cdot n$ changed to
$\omega^o \cdot (n-1)$.
\begin{quotation}
We started with a given constant function, which we characterized by a
determinate index.  We replaced a $\phi_b(\frak{z}, \frak{c}(b))$
within that functional by another functional, where the $\phi$-rank
decreased and the rank with respect to functions to the right of
$\phi$ [i.e., which come after $\phi$ in the order of definition] did
not increase.  Now we apply the same operation to the resulting
functional.  After finitely many steps we obtain a functional which
contains no function symbols at all, i.e., it is a numeral.

We have thus shown: a constant functional, which does not contain
$\epsilon$ and $\pi$, can be reduced to a numeral in finitely many
steps.\note{\pciteN{Ackermann:24}{18}}
\end{quotation}

\subsection{Ordinals, transfinite recursion, and finitism}

It is quite remarkable that the earliest extensive and detailed
technical contribution to the finitist project would make use of
transfinite induction in a way not dissimilar to Gentzen's later proof
by induction up to $\epsilon_0$.  This bears on a number of questions
regarding Hilbert's understanding of the strength of finitism.  In
particular, it is often said that Gentzen's proof is not finitist,
because it uses transfinite induction.  However, Ackermann's original
consistency proof for 2PRA$^-$ also uses transfinite induction, using an
index system which is essentially an ordinal notation system, just
like Gentzen's.  If it is granted that Ackermann's proof is
finitistic, but Gentzen's is not, i.e., transfinite induction up to
$\omega^{\omega^\omega}$ is finitistic but not up to $\epsilon_0$,
then where---and why---should the line be drawn?  Furthermore, the
consistency proof of 2PRA$^-$ is in essence a---putatively
finitistic---explanation of how to compute second order primitive
recursive functions, and a proof that the computation procedure
defined by them always terminates. In other words, it is a finitistic
proof that second order primitive recursive functions are well
defined.\note{\label{note-tait}\citeN{Tait:81} argues that finitism
coincides with primitive recursive arithmetic, and that therefore the
Ackermann function is not finitistic.  Tait does not present this as a
historical thesis, and his conceptual analysis remains unaffected by
the piece of historical evidence presented here. For further evidence
(dating however mostly from after 1931) see \pciteN{Zach:98}{\S5} and
Tait's response in \shortcite{Tait:00}.}
  
Ackermann was completely aware of the involvement of transfinite
induction in this case, but he sees in it no violation of the finitist
standpoint.
\begin{quotation}
The disassembling of functionals by reduction does not occur in the
sense that a finite ordinal is decreased each time an outermost
function symbol is eliminated.  Rather, to each functional corresponds
as it were a transfinite ordinal number as its rank, and the theorem,
that a constant functional is reduced to a numeral after carrying out
finitely many operations, corresponds to the other [theorem], that if
one descends from a transfinite ordinal number to ever smaller ordinal
numbers, one has to reach zero after a finite number of steps.
Now there is naturally no mention of transfinite sets or ordinal numbers
in our metamathematical investigations.  It is however interesting, that
the mentioned theorem about transfinite ordinals can be formulated so
that there is nothing transfinite about it.

Consider a transfinite ordinal number less than $\omega\cdot \omega$.
Each such ordinal number can be written in the form: $\omega\cdot
\frak{n} + \frak{m}$, where $\frak{n}$ and $\frak{m}$ are finite
numbers.  Hence such an ordinal can also be characterized by a pair of
finite numbers $(\frak{n}, \frak{m})$, where the order of these numbers
is of course significant.  To the descent in the series of ordinals
corresponds the following operation on the number pair $(\frak{n},
\frak{m})$.  Either the first number $\frak{n}$ remains the same, then
the number $\frak{m}$ is replaced by a smaller number $\frak{m}'$.  Or
the first number $\frak{n}$ is made smaller; then I can put an
arbitrary number in the second position, which can also be larger
than~$\frak{m}$.  It is clear that one has to reach the number pair
$(0, 0)$ after finitely many steps.  For after at most $\frak{m}+1$
steps I reach a number pair, where the first number is smaller than
$\frak{n}$.  Let $(\frak{n}', \frak{m}')$ be that pair.  After at most
$\frak{m}'+1$ steps I reach a number pair in which the first number is
again smaller than $\frak{n}'$, etc.  After finitely many steps one
reaches the number pair $(0,0)$ in this fashion, which corresponds to
the ordinal number~0.  In this form, the mentioned theorem contains
nothing transfinite whatsoever; only considerations which are
acceptable in metamathematics are used.  The same holds true if one
does not use pairs but triples, quadruples, etc.  This idea is not only
used in the following proof that the reduction of functionals
terminates, but will also be used again and again later on, especially
in the finiteness proof at the end of the
work.\note{\pciteN{Ackermann:24}{13--14}.}
\end{quotation}

Over ten years later, Ackermann discusses the application of
transfinite induction for consistency proofs in correspondence with
Bernays. Gentzen's consistency proof had been published
\cite{Gentzen:36}, and Gentzen asks, through Bernays,
\begin{quotation}
whether you [Ackermann] think that the method of proving termination
[\emph{Endlichkeitsbeweis}] by transfinite induction can be applied to
the consistency proof of your dissertation.  I would like it very
much, if that were possible.\note{``\sg [Gentzen fragt,] ob Sie der
Meinung sind, dass sich die Methode des Endlichkeitsbeweises durch
transfinite Induktion auf den Wf-Beweis Ihrer Dissertation anwenden
lasse. Ich w"urde es sehr begr"ussen, wenn das ginge.\eg'' Bernays to
Ackermann, November 27, 1936, Bernays Papers, ETH Z\"urich
Library/WHS, Hs 975.100.}
\end{quotation}
In his reply, Ackermann recalls his own use of transfinite ordinals in
the 1924 dissertation.
\begin{quotation}
I just realized now, as I am looking at my dissertation, that I
operate with transfinite ordinals in a similar fashion as
Gentzen.\note{``\sg Mir f"allt "ubrigens jetzt, wo ich gerade meine
Dissertation zur Hand nehme, auf, dass dort in ganz "ahnlicher Weise
mit transfiniten Ordnungszahlen operiert wird wie bei Gentzen.\eg''
Ackermann to Bernays, December 5, 1936, Bernays Papers, ETH Z\"urich
Library/WHS, Hs 975.101.}
\end{quotation}
A year and a half later, Ackermann mentions the transfinite induction
used in his dissertation again:
\begin{quotation}
I do not know, by the way, whether you are aware (I did at the time
not consider it as a transgression beyond the narrower finite
standpoint), that transfinite inferences are used in my dissertation.
(Cf., e.g., the remarks in the last paragraph on page 13 and the
following paragraph of my dissertation).\note{``\gq{Ich weiss "ubrigens
nicht, ob Ihnen bekannt ist (ich hatte das seiner Zeit nicht als
Ueberschreitung des engeren finiten Standpunktes empfunden), dass in
meiner Dissertation transfinite Schl"usse benutzt
werden. (Vgl. z.B. die Bemerkungen letzter Abschnitt Seite~13 und im
n"achstfolgenden Abschnitt meiner Dissertation.}'' Ackermann to
Bernays, June 29, 1938, Bernays papers, ETH-Z\"urich, Hs 975.114. The
passage Ackermann refers to is the one quoted above.}
\end{quotation}
These remarks may be puzzling, since they seem to suggest that Bernays
was not familiar with Ackermann's work.  This is clearly not the case.
Bernays corresponded with Ackerman extensively in the mid-20s about
the $\epsilon$-substitution method and the decision problem, and had
clearly studied Ackermann's dissertation.  Neither Bernays nor Hilbert
are on record objecting to the methods used in Ackermann's
dissertation.  It can thus be concluded that Ackermann's use of
transfinite induction was considered acceptable from the finitist
standpoint.

\subsection{The $\epsilon$-substitution method}\label{sec-esub}

As we have seen above, Hilbert had outlined an idea for a consistency
proof for systems involving $\epsilon$-terms already in early 1922
\cite{Hilbert:21a}, and a little more precisely in his talk of 1922
\cite{Hilbert:23} and in the 1922--23 lectures \cite{Hilbert:22c}.
Let us review the \emph{Ansatz} in the notation used in 1924: Suppose
a proof involves only one $\epsilon$ term $\epsilon_a A(a)$ and
corresponding \emph{critical formulas}
\[
A(\frak{k}_i) \to A(\epsilon_a A(a)),
\]
i.e., substitution instances of the transfinite axiom
\[
A(a) \to A(\epsilon_a A(a)).
\]
We replace $\epsilon_a A(a)$ everywhere with $0$, and transform the proof
as before by rewriting it in tree form (``dissolution into proof threads''),
eliminating free variables and evaluating numerical terms involving
primitive recursive functions.  Then the critical formulas take the form
\[
A(\frak{z}_i) \to A(0),
\]
where $\frak{z}_i$ is the numerical term to which $\frak{k}_i$
reduces.  A critical formula can now only be false if $A(\frak{z}_i)$
is true and $A(0)$ is false.  If that is the case, repeat the
procedure, now substituting $\frak{z}_i$ for $\epsilon_a A(a)$.  This
yields a proof in which all initial formulas are correct and no
$\epsilon$ terms occur.  

If critical formulas of the second kind, i.e.,
substitution instances of the induction axiom,
\[
\epsilon_a A(a) \neq 0 \to \ol{A(\delta \epsilon_a A(a))},
\]
also appear in the proof, the witness $\frak{z}$ has to be replaced
with the least $\frak{z}'$ so that $A(\frak{z}')$ is true.

The challenge was to extend this procedure to (a) cover more than one
$\epsilon$-term in the proof, (b) take care of nested
$\epsilon$-terms, and lastly (c) extend it to second-order
$\epsilon$'s and terms involving them, i.e, $\epsilon_f
\frak{A}_a(f(a))$.  This is what Ackermann set out to do in the last
part of his dissertation, and what he and Hilbert thought he had
accomplished.\note{The $\epsilon$-substitution method was subsequently
refined by \citeN{Neumann:27} and
\citeN{HilbertBernays:39}. \citeN{Ackermann:40} gives a consistency
proof for first-order arithmetic, using ideas of \citeN{Gentzen:36}.
See also \cite{Tait:65b} and \cite{Mints:94}.  Useful introductions to
the $\epsilon$-substitution method of \citeN{Ackermann:40} and to the
$\epsilon$-notation in general can be found in \citeN{Moser:00} and
\citeN{Leisenring:69}, respectively.}

The system for which Ackermann attempted to give a consistency proof
consisted of the system of second-order primitive recursive arithemtic
(see Section~\ref{sec-2pra} above) together with the transfinite
axioms:
\begin{center}
\begin{tabular}{lll}
1. & $A(a) \to A(\epsilon_a A(a))$ & $A_a f(a) \to A_a((\epsilon_f
A_b(f(b)))(a))$ \\ 2. & $A(\epsilon_a A(a)) \to \pi_a A(a) = 0$ &
$A_a(\epsilon_f A_b(f(b))(a)) \to \pi_f A_a(f(a)) = 0$\\ 3. &
$\ol{A(\epsilon_a A(a))} \to \pi_a A(a) = 1$ & $\ol{A_a(\epsilon_f
A_b(f(b))(a))} \to \pi_f A_a(f(a)) = 1$\\ 4. & $\epsilon_a A(a) \neq 0
\to \ol{A(\delta(\epsilon_a
A(a)))}$\note{\pciteN{Ackermann:24}{8}. The $\pi$-functions were
present in the earliest presentations in \cite{Hilbert:21a} as the
$\tau$-function and also occur in \cite{Hilbert:22c}.  They were
dropped from later presentations.}
\end{tabular}
\end{center}
The intuitive interpretation of $\epsilon$ and $\pi$, based on these
axioms, is obvious: $\epsilon_a \frak{A}(a)$ is a witness for
$\frak{A}(a)$ if one exists, and $\pi_a \frak{A}(a) = 1$ if
$\frak{A}(a)$ is false for all $a$, and $= 0$ otherwise.  The $\pi$
functions are not necessary for the development of mathematics in the
axiom system. They do, however, serve a function in the consistency
proof, viz., to keep track of whether a value of $0$ for $\epsilon_a
\frak{A}(a)$ is a ``default value'' (i.e., a trial substitution for
which $\frak{A}(a)$ may or may not be true) or an actual witness (a
value for which $\frak{A}(a)$ has been found to be true).

I shall now attempt to give an outline of the $\epsilon$-substitution
procedure defined by Ackermann. For simplicity, I will leave the case
of second-order $\epsilon$-terms (i.e., those involving $\epsilon_f$)
to the side.

An \emph{$\epsilon$-term} is an expression of the form $\epsilon_a
\frak{A}(a)$, where $a$ is the only free variable in~$\frak{A}$, and
similarly for a $\pi$-term.  For the purposes of the discussion below,
we will not specifically refer to $\pi$'s unless necessary, and most
definitions and operations apply equally to $\epsilon$-terms and
$\pi$-terms.  If a formula $A(a)$ or an $\epsilon$-term $\epsilon_a
\frak{A}(a)$ contains no variable-free subterms which are not
numerals, we call them \emph{canonical}.  Canonical formulas and
$\epsilon$-terms are indicated by a tilde: $\epsilon_a
\tilde{\frak{A}}(a)$.

The main notion in Ackermann's proof is that of a \emph{total
substitution}~S (\emph{Gesamtersetzung}).  It is a mapping of
canonical $\epsilon$- and $\pi$-terms to numerals and 0 or 1,
respectively. When canonical $\epsilon$-terms in a proof are
successively replaced by their values under the mapping, a total
substitution reduces the proof to one not containing any $\epsilon$'s.
If $S$ maps $\epsilon_a \tilde{\frak{A}}(a)$ to $\frak{z}$ and $\pi_a
\tilde{\frak{A}}(a)$ to $i$, then we say that $\tilde{\frak{A}}(a)$
receives a $(\frak{z}, \frak{i})$ substitution under~$S$ and write
$S(\tilde{\frak{A}}(a)) = (\frak{z}, \frak{i})$.  

It is of course not enough to define a mapping from the canonical
$\epsilon$-terms \emph{occuring in the proof} to numerals: The proof
may contain, e.g., $\epsilon_a \frak{A}(a, \phi(\epsilon_b
\frak{B}(b)))$.  To reduce this to a numeral, we first need a value
$\frak{z}$ for the term $\epsilon_b \frak{B}(b)$.  Replacing
$\epsilon_b \frak{B}(b)$ by $\frak{z}$, we obtain $\epsilon_a
\frak{A}(a, \phi(\frak{z}))$.  Suppose the value $\phi(\frak{z})$ is
$\frak{z}'$.  The total substitution then also has to specify a
substitution for $\epsilon_a \frak{A}(a, \frak{z}')$.

Given a total substitution~$S$, a proof is reduced to an
$\epsilon$-free proof as follows: First all $\epsilon$-free terms are
evaluated. (Such terms contain only numerals and primitive recursive
functions; these are computed and the term replaced by the numeral
corresponding to the value of the term) Now let $\frak{e}_{11}$,
$\frak{e}_{12}$, \dots{} be all the innermost (canonical) $\epsilon$-
or $\pi$-terms in the proof, i.e., all $\epsilon$- or $\pi$-terms
which do not themselves contain nested $\epsilon$- or $\pi$-terms
or constant (variable-fee) subterms which are not numerals.  The total
substitution specifies a numeral substitution for each of these.
Replace each $\frak{e}_{1i}$ by its corresponding numeral. Repeat the
procedure until the only remaining terms are numerals.  We write
$|\frak{e}|_S$ for the result of applying this procedure to the
expression (formula or term)~$\frak{e}$. Note that $|\frak{e}|_S$ is
canonical.

Based on this reduction procedure, Ackermann defines a notion of
subordination of canonical formulas.  Roughly, a formula
$\tilde{\frak{B}}(b)$ is subordinate to $\tilde{\frak{A}}(a)$ if in
the process of reducing some formula $\tilde{\frak{A}}(\frak{z})$, an
$\epsilon$-term $\epsilon_b \tilde{\frak{B}}(b)$ is replaced by a
numeral. For instance, $a = b$ is subordinate to $a = \epsilon_b(a =
b)$. Indeed, if $\tilde{\frak{A}}(a)$ is $a = \epsilon_b (a = b)$,
then the reduction of $\tilde{\frak{A}}(\frak{z}) \equiv \frak{z} =
\epsilon_b (\frak{z} = b)$ would use a replacement for the
$\epsilon$-term belonging to $\tilde{\frak{B}}(\frak{z} = b)$.\note{It
is not clear whether the definition is supposed to apply to the
formulas with free variables (i.e., to $a = b$ and $a = \epsilon_b(a =
b)$ in the example) or to the corresponding substituion instances.
The proof following the definition on p.~21 of \cite{Ackermann:24}
suggests the former, however, later in the procedure for defining a
sequence of total substitutions it is suggested that the
$\epsilon$-expressions corresponding to formulas subordinate to
$\tilde{\frak{A}}(a)$ receive substitutions---but according to the
definition of a total substitution only $\epsilon$-\emph{terms}
($\epsilon_b (\frak{z} = b)$ in the example) receive substitutions.}
It is easy to see that this definition corresponds to the notion of
subordination as defined in \citeN{HilbertBernays:39}.  An
$\epsilon$-\emph{expression} is an expression of the form $\epsilon_a
\frak{A}(a)$. If $\epsilon_a \frak{A}(a)$ contains no free variables,
it is called an $\epsilon$-\emph{term}.  If an $\epsilon$-term
$\epsilon_b \frak{B}(b)$ occurs in an expression (and is different
from it), it is said to be \emph{nested} in it.  If an
$\epsilon$-expression $\epsilon_b \frak{B}(a, b)$ occurs in an
expression in the scope of $\epsilon_a$, then it is \emph{subordinate}
to that expression.  Accordingly, we can define the degree of an
$\epsilon$-term and the rank of an $\epsilon$-expression as follows:
An $\epsilon$-term with no nested $\epsilon$-subterms is of degree~1;
otherwise its degree is the maximum of the degrees of its nested
$\epsilon$-subterms $+ 1$. The rank of an $\epsilon$-expression with
no subordinate $\epsilon$-expressions is~1; otherwise it is the
maximum of the ranks of its subordinate $\epsilon$-expressions~$+1$.
If $\tilde{\frak{B}}(b)$ is subordinate to $\tilde{\frak{A}}(a)$
according to Ackermann's definition, then $\epsilon_b
\tilde{\frak{B}}(b)$ is subordinate in the usual sense to $\epsilon_a
\tilde{\frak{A}}(a)$, and the rank of $\epsilon_b \tilde{\frak{B}}(b)$
is less than that of $\epsilon_a \tilde{\frak{A}}(a)$.  The notion of
degree corresponds to an ordering of canonical formulas used for the
reduction according to a total substitution in Ackermann's procedural
definition: First all $\epsilon$-terms of degree~1 (i.e., all
innermost $\epsilon$-terms) are replaced, resulting (after evaluation
of primitive recursive functions) in a partially reduced proof.  The
formulas corresponding to innermost $\epsilon$-terms now are reducts
of $\epsilon$-terms of degree~2 in the original proof. The canonical
formulas corresponding to $\epsilon$-terms of degree~1 are called the
formulas of \emph{level}~1, the canonical formulas corresponing to the
innermost $\epsilon$-terms in the results of the first reduction step
are the formulas of level~2, and so forth.

The consistency proof proceeds by constructing a sequence $S_1$,
$S_2$, \dots{} of total substitutions together with bookkeeping
functions $f_i(\tfA(a), j) \to \{0, 1\}$,\note{The bookkeeping
functions are introduced here and are not used by Ackermann.  The
basic idea is that that in case~(3), substitutions for some formulas
are discarded, and the next substitution is given the ``last'' total
substitution where the substitution for the formula was not yet marked
as discarded. Instead of explicit bookkeeping, Ackermann uses the
notion of a formula being ``remembered'' as having its value not
discarded.} which eventually results in a \emph{solving substitution},
i.e., a total substitution which reduces the proof to one which
contains only correct $\epsilon$-free formulas.  We begin with a total
substitution~$S_1$ which assigns $(0, 1)$ to all canonical formulas,
and set $f_1(\tfA(a), 1) = 1$ for all $\tfA(a)$ for which $S_1$
assigns a value.  If $S_i$ is a solving substitution, the procedure
terminates. Otherwise, the next total substitution~$S_{i+1}$ is
obtained as follows: If $S_i$ is not a solving substitution, at least
one of the critical formulas in the proof reduces to an incorrect
formula.  We have three cases:

\begin{enumerate}
\item Either an $\epsilon$-axiom $\frak{A}(a) \to \frak{A}(\epsilon_a
    \frak{A}(a))$ or a $\pi$-axiom of the first kind
    $\frak{A}(\epsilon_a \frak{A}(a)) \to \pi_a \frak{A}(a) = 0$
    reduces to a false formula of the form $\tilde{\frak{A}}(\frak{z})
    \to \tilde{\frak{A}}(0)$ or $\tilde{\frak{A}}(0) \to 1 = 0$, and
    $S_i(\tilde{\frak{A}}(a)) = (0, 1)$.  Pick one such
    $\tilde{\frak{A}}(a)$ of lowest level (i.e., $\epsilon_a
    \frak{A}(a)$ of lowest degree).  

    If $\tfA(0) \to 1 = 0$ is incorrect, $\tilde{\frak{A}}(0)$ must be
    correct; let $S_{i+1}(\tilde{A}(a)) = (0, 0)$. Otherwise
    $\tfA(\fz) \to \tfA(0)$ is incorrect and hence $\tfA(\fz)$ must be
    correct; then let $S_{i+1}(\tfA(a)) = (\fz, 0)$.  In either case,
    set $f_{i+1}(\tilde{\frak{A}}(a), i+1) = 1$.
    
    For other formulas $\tfB(b)$, $S_{i+1}(\tfB(b)) =
    S_{j}(\tilde{\frak{B}}(b))$ where $j$ is the greatest index $\le
    i$ such that $f_i(\tfB(b), j) = 1$. $S_{i+1}(\tfB(b)) = (0, 1)$ if
    no such $j$ exists (i.e., $\tfB(b)$ has never before received an
    example substitution). Also, let $f_{i+1}(\tfB(b), i+1) = 1$.  For
    all canonical formulas $\tilde{\frak{C}}(c)$, let
    $f_{i+1}(\tilde{\frak{C}}(c), j) = f_{i}(\tilde{\frak{C}}(c), j)$
    for $j \le i$.

\item Case (1) does not apply, but at least one of the minimality
axioms $\epsilon_a \frak{A}(a) \neq 0 \to
\ol{\frak{A}}(\delta(\epsilon_a \frak{A}(a))$ reduces to a false
formula, $\frak{z} \neq 0 \to \ol{\tilde{\frak{A}}}(\frak{z} - 1)$.
This is only possible if $S_i(\tilde{\frak{A}}(a)) = (\frak{z}, 0)$.
Again, pick the one of lowest level, and let
$S_{i+1}(\tilde{\frak{A}}(a)) = (\frak{z}-1, 0)$ and $f_{i+1}(\tfA(a),
i+1) = 1$.  Substitutions for other formulas and bookkeeping functions
are defined as in case~(1).
    
\item Neither case (1) nor (2) applies, but some instance of an
$\epsilon$-axiom of the form $\frak{A}(a) \to \frak{A}(\epsilon_a
\frak{A}(a))$ or of a $\pi$-axiom of the form
$\ol{\frak{A}}(\epsilon_a \frak{A}(a)) \to \pi_a \frak{A}(a) = 1$,
e.g., $\tilde{\frak{A}}(\frak{a}) \to \tilde{\frak{A}}(\frak{z}))$ or
$\ol{\tilde{\frak{A}}}(\frak{z}) \to 0 = 1$, reduces to an incorrect
formula. We then have $S_i(\tilde{\frak{A}}(a)) = (\frak{z}, 0)$
(since otherwise case (1) would apply).  In either case,
$|\tilde{\frak{A}}(\frak{z}))|_{S_i}$ must be incorrect.  Let $j$ be
the least index where $S_j(\tilde{\frak{A}}(a)) = (\frak{z}, 0)$ and
$f_i(\tilde{\frak{A}}(a)), j) = 1$.  At the preceding total
substitution $S_{j-1}$, $S_j(\tilde{\frak{A}}(a)) = (0, 1)$ or
$S_{j-1}(A(a)) = (\frak{z} + 1, 0)$, and
$|\frak{A}(\frak{z})|_{S_{j-1}}$ is correct.
$\tilde{\frak{A}}(\frak{z})$ thus must reduce to different formulas
under $S_{j-1}$ and under $S_i$, which is only possible if a formula
subordinate to $\tilde{\frak{A}}$ reduces differently under $S_{j-1}$
and $S_i$.  

For example, suppose $\tilde{\frak{A}}(a)$ is really
$\tilde{\frak{A}}(a, \epsilon_b \tilde{\frak{B}}(a, b))$. Then the
corresponding $\epsilon$-axiom would be
\[
\tfA(a,\epsilon_b \tfB(a, b)) \to \tfA(\underbrace{\epsilon_a \tfA(a,
\epsilon_b \tfB(a, b))}_{\epsilon_a \tfA(a)}, \epsilon_b
\tfB(\underbrace{\epsilon_a \tfA(a, \epsilon_b\tfB(a,b))}_{\epsilon_a
\tfA(a)}, b))
\]
An instance thereof would be
\[
\tfA(\frak{a},\epsilon_b \tfB(\frak{a}, b) \to \tfA(\epsilon_a \tfA(a,
\epsilon_b \tfB(a, b)), \epsilon_b \tfB(\epsilon_a \tfA(a,
\epsilon_b\tfB(a,b)), b)).
\] 
This formula, under a total substitution with $S_i(\tfA(a, \epsilon_b \tfB(a,
b))) = (\frak{z}, 0)$ reduces to
\[
\tfA(\frak{a},\epsilon_b \tfB(\frak{a}, b) \to \tfA(\frak{z},
\epsilon_b \tfB(\frak{z}, b))
\] 
The consequent of this conditional, i.e., $\tfA(\frak{z})$,
can reduce to different formulas under $S_i$ and $S_{j-1}$ only if
$\epsilon_b \tfB(\frak{z}, b)$ receives different substitutions under
$S_i$ and $S_{j-1}$, and $\tfB(a, b)$ is subordinate to $\tfA(a)$.
 
    The next substitution is now defined as follows: Pick an innermost
    formula subordinate to $\tfA(a)$ which changes substitutions, say
    $\tilde{\frak{B}}(b)$.  For all formulas $\tilde{\frak{C}}(c)$
    which are subordinate to $\tilde{\frak{B}}(b)$ as well as
    $\tilde{\frak{B}}(b)$ itself, we set $f_{i+1}(\tilde{\frak{C}}(c),
    k) = 1$ for $j \le k \le i+1$ and $f_{i+1}(\tilde{\frak{C}}(c), k)
    = 0$ for all other formulas.  For $k < j$ we set
    $f_{i+1}(\tilde{\frak{C}}(c), n) = f_i(\tilde{\frak{C}}(c), k)$
    for all $\tilde{\frak{C}}(c)$.  The next substitution $S_{i+1}$ is
    now given by $S_{i+1}(\tilde{\frak{C}}(c)) =
    S_k(\tilde{\frak{C}}(c))$ for $k$ greatest such that
    $f_{i+1}(\tilde{\frak{C}}(c), k) = 1$ or $= (0, 1)$ if no such $k$
    exists.
\end{enumerate}
Readers familiar with the substitution method defined
in \citeN{Ackermann:40} will note the following differences:
\begin{enumerate}
\item[a.] \citeN{Ackermann:40} uses the notion of a \emph{type} of an
$\epsilon$-term and instead of defining total substitutions in terms
of numeral substitutions for canonical $\epsilon$-terms, assigns a
function of finite support to $\epsilon$-types.  This change is merely
a notational convenience, as these functional substitutions can be
recovered from the numeral substitutions for canonical
$\epsilon$-terms.  For example, if $S$ assigns the substitutions to the
canonical terms on the left, then a total substitution in the sense of
\citeN{Ackermann:40} would assign the function $g$ on the right to the
type $\epsilon_a \tfA(a, b)$:
\[
\begin{array}{l@{\qquad}l}
S(\tfA(a, \fz_1) = \fz_1' & g(\epsilon_a \tfA(a, b))(\fz_1) = \fz_1' \\
S(\tfA(a, \fz_3) = \fz_2' & g(\epsilon_a \tfA(a, b))(\fz_2) = \fz_2' \\
S(\tfA(a, \fz_2) = \fz_3' & g(\epsilon_a \tfA(a, b))(\fz_3) = \fz_3' \\
\end{array}
\]
\item[b.] In case (2), dealing with the least number (induction)
axiom, the next substitution is defined by reducing the substituted
numeral~$\fz$ by~1, whereas in \cite{Ackermann:40}, we immediately
proceed to the least~$\fz'$ such that $\tfA(\fz')$ is correct.  This
makes the procedure converge more slowly, but also suggests that in
certain cases (depending on which other critical formulas occur in the
proof), the solving substitution does not necessarily provide example
substitutions which are, in fact, least witnesses.  
\item[c.] The main difference in the method lies in case~(3). Whereas
in \shortcite{Ackermann:40}, example substitutions for all
$\epsilon$-types of rank lower than that of the changed $\epsilon_a
\tfA(a)$ are retained, and all others are reset to initial
substitutions (functions constant equal to $0$), in
\shortcite{Ackermann:24}, only the substitutions of some
$\epsilon$-terms actually subordinate to $\epsilon_a \tfA(a)$ are
retained, while others are not reset to initial substitutions,
but to substitutions defined at some previous stage.
\end{enumerate}

\subsection{Assessment and complications}

A detailed analysis of the method and of the termination proof given
in the last part of Ackermann's dissertation must wait for another
occasion, if only for lack of space.  A preliminary assessment can,
however, already be made on the basis of the outline of the
substitution process above.  Modulo some needed clarification in the
definitions, the process is well-defined and terminates at least for
proofs containing only least-number axioms (critical formulas
corresponding to axiom~(4)) of rank~1.  The proof that the procedure
terminates (\S9 of \citeN{Ackermann:24}) is opaque, especially in
comparison to the proof by transfinite induction for primitive
recursive arithmetic.  The definition of a substitution method for
second-order $\epsilon$-terms is insufficient, and in hindsight it is
clear that a correct termination proof for this part could not have
been given with the methods available.\note{With the restriction on
second-order $\epsilon$-terms imposed by Ackermann, and discussed
below, the system for which a consistency proof was claimed is
essentially elementary analysis, a predicative system. A consistency
proof using the $\epsilon$-substitution method for this system was
given by \citeN{Mints:96}.}

Leaving aside, for the time being, the issue of what was
\emph{actually} proved in \citeN{Ackermann:24}, the question remains
of what was \emph{believed} to have been proved at the time.  The
system, as given in \shortcite{Ackermann:24}, had two major
shortcomings: A footnote, added in proof, states:
\begin{quotation}
[The formation of $\epsilon$-terms] is restricted in that a function
variable $f(a)$ may not be substituted by a functional $\frak{a}(a)$,
in which $a$ occurs in the scope of an
$\epsilon_f$.\note{\pciteN{Ackermann:24}{9}}
\end{quotation}
This applies in particular to the second-order $\epsilon$-axioms
\[A_a(f(a)) \to A_a(\epsilon_f A_b(f(b))(a)).\] If we view $\epsilon_f
A_b(f(b))(a)$ as the function ``defined by'' $A$, and hence the
$\epsilon$-axiom as the $\epsilon$-calculus analog of the
comprehension axiom, this amounts roughly to a restriction to
arithmetic comprehension, and thus a predicative system. This
shortcoming, and the fact that the restriction turns the system into a
system of predicative mathematics was pointed out by
\citeN{Neumann:27}.

A second lacunae was the omission of an axiom of
$\epsilon$-extensionality for second-order $\epsilon$-terms, i.e.,
\[
(\forall f)(A(f) \rightleftharpoons B(f)) \to \epsilon_f A(f) =
\epsilon_f B(f),
\]
which corresponds to the axiom of choice.  Both problems were the
subject of correspondence with Bernays in 1925.\note{Ackermann to
Bernays, June~25, 1925, Bernays Papers, ETH Z\"urich, Hs. 975.96.} A
year later, Ackermann is still trying to extend and correct the proof,
now using $\epsilon$-types:
\begin{quotation}
I am currently working again on the $\epsilon_f$-proof and am pushing
hard to finish it.  I have already told you that the problem can be
reduced to one of number theory.  To prove the number-theoretic
theorem seems to me, however, equally hard as the problem itself.  I am
now again taking the approach, which I have tried several times
previously, to extend the definition of a ground type so that even
$\epsilon$ with free function variables receive a substitution.  This
approach seems to me the most natural, and the equality axioms
$(f)(A(f) \rightleftharpoons B(f)) \to \epsilon_f Af \equiv \epsilon_f
Bf$ would be treated simultaneously. I am hopeful that the obstacles
previously encountered with this method can be avoided, if I use the
$\epsilon_a$ formalism and use substitutions for the $\epsilon_f$
which may contain $\epsilon_a$ instead of functions defined without
$\epsilon$.  I have, however, only thought through some simple special
cases.\note{``\gq{Ich habe augenblicklich den $\epsilon_f$-Beweis
wieder vorgenommen, und versuche mit aller Gewalt da zum Abschlu"s zu
kommen. Da"s sich das Problem auf ein zahlentheoretisches reduzieren
l"a"st, hatte ich Ihnen damals ja schon mitgeteilt. Den
zahlentheoretischen Satz allgemein zu beweisen scheint mir aber ebenso
schwierig wie das ganze Problem. Ich habe nun den schon mehrfach von
mir versuchten Weg wieder eingeschlagen, den Begriff des Grundtyps so
zu erweitern, das auch die $\epsilon$ mit freien Funktionsvariablen
eine Ersetzung bekommen.  Dieser Weg scheint ja auch der
nat"urlichste, und die Gleichheitsaxiome $(f)(A(f) \rightleftharpoons
B(f)) \to \epsilon_f Af \equiv \epsilon_f Bf$ w"urden dann gleich
mitbehandelt. Ich habe einige Hoffnung, da"s die sich fr"uher auf
diesem Weg einstellenden Schwierigkeiten vermieden werden k"onnen,
wenn ich den $\epsilon_a$-Formalismus benutze und statt ohne
$\epsilon$ definierte Funktionen, solche zur Ersetzung f"ur die
$\epsilon_f$ nehme, die ein $\epsilon_a$ enthalten k"onnen. Ich habe
mir aber erst einfache Spezialf"alle "uberlegt.}'' Ackermann to
Bernays, March 31, 1926. ETH Z"urich/WHS, Hs~975.97. Although
Ackermann's mention of ``ground types'' precedes the publication of
\citeN{Neumann:27}, the latter paper was submitted for publication
alread on July 29, 1925.}
\end{quotation}
In 1927, Ackermann developed a second proof of
$\epsilon$-substitution, using some of von Neumann's ideas (in
particular, the notion of an $\epsilon$-type, \emph{Grundtyp}). The
proof is unfortunately not preserved in its entirety, but references
to it can be found in the correspondence.  On April 12, 1927, Bernays
writes to Ackermann:
\begin{quotation}
Finally I have thought through your newer proof for consistency of the
$\epsilon_a$'s based on what you have written down for me befor your
departure, and believe that I have seen the proof to be
correct.\note{``\gq{Letzthin habe ich mir Ihren neueren Beweis der
Widerspruchsfr[eiheit] f"ur die $\epsilon_a$ an Hand dessen, was Sie
mir vor Ihrer Abreise aufschrieben, genauer "uberlegt und glaube
diesen Beweis als richtig eingesehen zu haben.}'' Bernays to
Ackermann, April 12, 1927, in the posession of Hans Richard Ackermann.
Bernays continues to remark on specifics of the proof, roughly, that
when example substitutions for $\epsilon$-types are revised (the
situation corresponding to case (3) in Ackermann's original proof),
the substitutions for types of higher rank have to be reset to the
initial substitution.  He gives an example that shows that if this is
not done, the procedure does not terminate. He also suggests that it
would be more elegant to treat all types of the same rank at the same
time and gives an improved estimate for the number of steps necessary.
Note that the reference to ``$\epsilon_a$'s'' (as opposed to
$\epsilon_f$) suggest that the proof was only for the first-order
case.  A brief sketch of the proof is also contained in a letter from
Bernays to Weyl, dated January 5, 1928 (ETH Z"urich/WHS, Hs.~91.10a).}
\end{quotation}
Ackermann also refers to the proof in a letter to Hilbert from 1933:
\begin{quotation}
As you may recall, I had at the time a second proof for the
consistency of the $\epsilon_a$'s.  I never published that proof, but
communicated it to Prof.~Bernays orally, who then verified it.
Prof.~Bernays wrote to me last year that the result does not seem to
harmonize with the work of G"odel.\note{``\gq{Wie Sie sich vielleicht
erinnern, hatte ich damals einen 2.~Beweis f"ur die
Widerspruchsfreiheit der $\epsilon_a$. Dieser Beweis ist von mir nie
publiziert worden, sondern nur Herrn Prof. Bernays m"undlich
mitgeteilt worden, der sich auch damals von seiner Richtigkeit
"uberzeugte.  Prof. Bernays schrieb mir nun im vergangenen Jahre, da"s
das Ergebnis ihm mit der G"odelschen Arbeit nicht zu harmonisieren
scheine.}'' Ackermann to Hilbert, August 23, 1933, Hilbert-Nachla"s,
Nieders"achsische Staats- und Universit"atsbibliothek,
Cod. Ms. Hilbert 1.  Ackermann did not then locate the difficulty, and
even a year and a half later (Ackermann to Bernays, December 8, 1934,
ETH Z"urich/WHS, Hs 975.98) suggested a way that a finitistic
consistency proof of arithmetic could be found based on work of
Herbrand and Bernays's drafts for the second volume of
\emph{Grundlagen}.}
\end{quotation}
In Hilbert's address to the International Congress of Mathematicians
in 1928 \cite{Hilbert:28}, the success of Ackermann's and von
Neumann's work on $\epsilon$-substitution for first-order systems is
also taken for granted. Although Hilbert poses the extension of the
proof to second-order systems as an open problem, there seems no doubt
in his mind that the solution is just around the
corner.\note{``Problem I. The consistency proof of the
$\epsilon$-axiom for the function variable~$f$. We have the outline of
a proof. Ackermann has already carried it out to the extent that the
only remaining task consists in the proof of an elementary finiteness
theorem that is purely arithmetical.''  \citeN{Hilbert:28}, translated
in \pciteN{Mancosu:99}{229}. The extension to
$\epsilon$-extensionality is Problem~III.}

It might be worthwhile to mention at this point that at roughly the
same time a third attempt to find a satisfactory consistency proof was
made.  This attempt was based not on $\epsilon$-substituion, but on
Hilbert's so-called unsuccessful proof (\emph{verungl"uckter Beweis}).
\begin{quotation}
While working on the \emph{Grundlagenbuch}, I found myself motivated
to re-think Hilbert's second consistency proof for the
$\epsilon$-axiom, the so-called ``unsuccessful'' proof, and it now
seems to me that it can be fixed after all.\note{``\gq{Anl"asslich der
Arbeit f"ur das Grundlagenbuch sah ich mich dazu angetrieben, den
zweiten Hilbertschen Wf.-Beweis f"ur das $\epsilon$-Axiom, den
sogenannten "`verungl"uckten"' Beweis, nochmals zu "uberlegen, und es
scheint mir jetzt, dass dieser sich doch richtig stellen l"asst.}''
Bernays to Ackermann, October 16, 1929, in the possession of Hans
Richard Ackermann.  Bernays continues with a detailed exposition of
the proof, but concludes that the proof probably cannot be extended to
include induction, for which $\epsilon$-substitution seems better
suited.}
\end{quotation}
This proof bears a striking resemblance to the proof the first
$\epsilon$-theorem in \cite{HilbertBernays:39} and to a seven-page
sketch in Bernays's hand of a ``consistency proof for the logical
axiom of choice'' found bound with lecture notes to Hilbert's course
on ``Elements and principles of mathematics'' of 1910.\note{The sketch
bears the title ``Wf.-Beweis f"ur das logische Auswahl-Axiom'', and is
inserted in the front of \emph{Elemente und Prinzipienfragen der
Mathematik}, Sommer-Semester 1910. Library of the Mathematisches
Institut, Universit"at G"ottingen, 16.206t14.  A note in Hilbert's
hand says ``Einlage in W.S. 1920.'' However, the $\epsilon$-Axiom used
is the more recent version $Ab \to A\epsilon_a A_a$ and not the
original, dual $A\epsilon_a Aa \to Ab$.  It is thus very likely that
the sketch dates from after 1923.}  This ``unsuccessful'' proof seems
to me to be another but independent contribution to the development of
logic and the $\epsilon$-calculus, independent of the substituion
method.  Note that Bernays's proof of Herbrand's theorem in
\cite{HilbertBernays:39} is based on the (second) $\epsilon$-theorem
is the first correct published proof of that important result.

The realization that the consistency proof even for first-order
$\epsilon$'s was problematic came only with G\"odel's incompleteness
results.  In a letter dated March~10, 1931, von Neumann presents an
example that shows that in the most recent version of Ackermann's
proof, the length of the substitution process not only depends on the
rank and degree of $\epsilon$-terms occurring in the proof, but also
on numerical values used as substitutions. He concludes:
\begin{quotation}
I think that this answers the question, which we recently discussed
when going through Ackermann's modified proof, namely whether an
estimate of the length of the correction process can be made uniformly
and independently of numerical substituends, in the negative.  At this
point the proof of termination of the procedure (for the next higher
degree, i.e., 3) has a gap.\note{``\gq{Ich glaube, dass damit die
Frage, die wir bei der Durchsprechung des modifizierten Ackermannschen
Beweises zuletzt diskutierten, ob n"amlich eine L"angen-Absch"atzung
f"ur das Korrigier-Verfahren unabh"angig von der Gr"osse der
Zahlen-Substituenden gleichm"assig m"oglich sei, verneinend
beantwortet ist. An diesem Punkte ist dann der Nachweis des endlichen
Abbrechens dieses Verfahrens (f"ur den n"achsten Grad, d.h. 3)
jedenfalls l"uckenhaft.)}'' von Neumann to Bernays, March 10, 1931,
Bernays Papers, ETH Z"urich/WHS, Hs.~975.3328. Von Neumann's example
can be found in \pciteN{HilbertBernays:39}{123}.}
\end{quotation}
There is no doubt that the discussion of the consistency proof was
precipitated by G\"odel's results, as both von Neumann and Bernays
were aware of these results, and at least von Neumann realized the
implications for Hilbert's Program and the prospects of a finitistic
consistency proof for arithmetic.  Bernays corresponded with G"odel on
the relevance of G"odel's result for the viability of the project of
consistency proofs just before and after von Neumann's counterexample
located the difficulty in Ackermann's proof.  On January 18, 1931,
Bernays writes to G"odel:
\begin{quotation}
If one, as does von Neumann, assumes as certain that any finite
consideration can be formulated in the framework of System~$P$---I
think, as you do too, that this is not at all obvious---one arrives at
the conclusion that a finite proof of consistency of $P$ is
impossible.
\end{quotation}
The puzzle, however, remained unresolved for Bernays even after von
Neumann's example, as he writes to G"odel just after the exchange with
von Neumann, on April 20, 1931:
\begin{quotation}
The confusion here is probably connected to that about Ackermann's
proof for the consistency of number theory (System~$\frak{Z}$), which
I have not so far been able to clarify.  

That proof---on which Hilbert has reported in his Hamburg talk on the
``foundations of mathematics''\note{\citeN{Hilbert:28}.} [\dots]---I
have repeatedly thought through and found correct.  On the basis of
your results one must now conclude that this proof cannot be
formalized within System~$\frak{Z}$; indeed, this must hold even if
one restricts the system whose consistency is to be proved by leaving
only addition and multiplication as recursive definitions. On the
other hand, I do not see which part of Ackermann's proof makes the
formalization within $\frak{Z}$ impossible, in particular if the
problem is so restricted.\note{In a letter dated May 3, 1931, Bernays
suggests that the problem lies with certain types of recursive
definitions The Bernays--G"odel correspondence will shortly be
published in Volume~IV of G"odel's collected works.  For more on the
reception of G\"odel's results by Bernays and von Neumann, see
\citeN{Dawson:88} and \citeN{Mancosu:99a}.}
\end{quotation}
G"odel's results thus led Bernays, and later Ackermann to reexamine
the methods used in the consistency proofs.  A completion of the
project had to wait until 1940, when Ackermann was able to carry
through the termination proof based on transfinite
induction---following \citeN{Gentzen:36}---on~$\epsilon_0$.

\section{Conclusion}

With the preceding exposition and analysis of the development of
axiomatizations of logic and mathematics and of Hilbert and
Ackermann's consistency proofs I hope to have answered some open
questions regarding the historical development of Hilbert's Program.
Hilbert's \emph{technical} project and its evolution is without doubt
of tremendous importance to the history of logic and the foundations
of mathematics in the 20th century.  Moreover, an understanding of the
technical developments can help to inform an understanding of the
history and prospects of the \emph{philosophical} project.  The
lessons drawn in the discussion, in particular, of Ackermann's use of
transfinite induction, raise more questions.  The fact that
transfinite induction in the form used by Ackermann was so readily
accepted as finitist, not just by Ackermann himself, but also by
Hilbert and Bernays leaves open two possibilities: either they were
simply wrong in taking the finitistic nature of Ackermann's proof for
granted and the use of transfinite induction simply cannot be
reconciled with the finitist standpoint as characterized by Hilbert
and Bernays in other writings, or the common view of what Hilbert
thought the finitist standpoint to consist in must be revised.
Specifically, it seems that the explanation of why transfinite
induction is acceptable stresses one aspect of finitism while
downplaying another: the \emph{objects} of finitist reasoning
are---finite and---intuitively given, whereas the methods of proof
were not required to have the epistemic strength that the finitist
standpoint is usually thought to require (i.e., to guarantee, in one
sense or another, the intuitive evidence of the resulting theorems).
Of course, the question of whether Hilbert can make good on his claims
that finitistic reasoning affords this intutive evidence of its
theorems is one of the main difficulties in a philosphical assessment
of the project (see, e.g., \citeN{Parsons:98}).

I have already hinted at the implications of a study of the practice
of finitism for philosophical reconstructions of the finitist view (in
note~\ref{note-tait}).  We are of course free to latch on to this or
that aspect of Hilbert's ideas (finitude, intuitive evidence, or
surveyability) and develop a philosophical view around it.  Such an
approach can be very fruitful, and have important and inightful
results (as, e.g, the example of Tait's~\shortcite{Tait:81} work
shows).  The question is to what extent such a view should be accepted
as a reconstruction of Hilbert's view as long as it makes the practice
of the technical project come out off base.  Surely rational
reconstruction is governed by something like a principle of charity.
Hilbert and his students, to the extent possible, should be construed
so that what they preached is reflected in their practice.  This
requires, of course, that we know what the practice was.  If nothing
else, I hope to have provided some of the necessary data for that.

\acknowledgements

Parts of this paper were presented to the Townsend Center Working
Group in History and Philosophy of Logic and Mathematics at the
University of Caliornia, Berkeley, December 1999, at a workshop on
Hilbert at the Institut d'Histoire des Sciences et Techniques of the
Universit\'e Paris~I (partially funded by the Townsend Center through
an International Working Group grant), May 2000, at he International
Symposium on the History of Logic at the University of Helsinki, June
2000, and at the Department of Logic and Philosophy of Science of the
University of California, Irvine, January 2001.  I am indebted to Hans
Richard Ackermann for providing me with copies of Wilhelm Ackermann's
correspondence, and to Wilfried Sieg for providing me with Kneser's
notes to the 1921--22 and 1922--23 lectures
\cite{Hilbert:21a,Hilbert:22c}.  I would also like to thank Grigori
Mints and Paolo Mancosu for valuable comments.  My understanding of
the $\epsilon$-substitution method in general and Ackermann's proof in
particular has benefitted greatly from conversations with and writings
of Grigori Mints, Georg Moser, and W. W. Tait. Any errors in my
presentation are, of course, mine.

\theendnotes

\end{article}
\end{document}